\documentclass[trsc,nonblindrev]{informs3} 

\OneAndAHalfSpacedXI 

\usepackage{adjustbox}
\usepackage{natbib}
\bibpunct[, ]{(}{)}{,}{a}{}{,}%
\TheoremsNumberedThrough     

\EquationsNumberedThrough    

\usepackage{makecell}
\usepackage{booktabs}
\usepackage[figuresright]{rotating}
\usepackage{hyperref}
\hypersetup{
    colorlinks,
    linkcolor={red!50!black},
    citecolor={blue!50!black},
    urlcolor={blue!80!black}
}
\usepackage{multirow}
\usepackage{algorithm}
\usepackage[noend]{algpseudocode}
\usepackage{xcolor}
\usepackage{lipsum}
\usepackage{amsfonts}

\begin{document}

    \RUNAUTHOR{Lai et al.}
    \RUNTITLE{PRP with Speed Optimization and Uneven Topography}
    \TITLE{The Pollution-Routing Problem with Speed Optimization and Uneven Topography}

    \ARTICLEAUTHORS{%
        \AUTHOR{David Lai}
        \AFF{School of Industrial Engineering, Eindhoven University of Technology, 5600 MB Eindhoven, Netherlands,\\ \EMAIL{s.w.lai@tue.nl}}
   \AUTHOR{Yasel Costa}
        \AFF{MIT-Zaragoza International Logistics Program, Zaragoza Logistics Center, Spain\\ \EMAIL{ycosta@zlc.edu.es}}
           \AUTHOR{Emrah Demir}
        \AFF{Cardiff Business School, Cardiff University, Cardiff, United Kingdom\\ \EMAIL{demire@cardiff.ac.uk}}
           \AUTHOR{Alexandre M. Florio, Tom Van Woensel}
      \AFF{School of Industrial Engineering, Eindhoven University of Technology, 5600 MB Eindhoven, Netherlands,\\ \EMAIL{a.de.macedo.florio@tue.nl}, \EMAIL{t.v.woensel@tue.nl}}
    } 

\ABSTRACT{%
This paper considers a joint pollution-routing and speed optimization problem (PRP-SO) where fuel costs and $\text{CO}_2e$ emissions depend on the vehicle speed, arc payloads, and road grades. We present two methods, one approximate and one exact, for solving the PRP-SO. The approximate strategy solves large-scale instances of the problem with a tabu search-based metaheuristic coupled with an efficient fixed-sequence speed optimization algorithm. The second strategy consists of a tailored branch-and-price (BP) algorithm in which speed optimization is managed within the pricing problem. We test both methods on modified Solomon benchmarks and newly constructed real-life instance sets. Our BP algorithm solves most instances with up to 50 customers and many instances with 75 and 100 customers. The heuristic is able to find near-optimal solutions to all instances and requires less than one minute of computational time per instance. Results on real-world instances suggest several managerial insights. First, fuel savings of up to 53\% are realized when explicitly taking into account arc payloads and road grades. Second, fuel savings and emissions reduction are also achieved by scheduling uphill customers later along the routes. Lastly, we show that ignoring elevation information when planning routes leads to highly inaccurate fuel consumption estimates.

}%


    \maketitle
\clearpage\newpage

\section{Introduction}
Road freight transportation is vital to the functioning of the economy and the supply chain. However, significant negative impacts on people and the environment due to excessive energy usage and considerable greenhouse emissions need to be considered. According to the International Energy Agency \citep{IEA:2020}, global transportation is still responsible for 24\% of direct $\text{CO}_2$-equivalent ($\text{CO}_2e$) emissions  from fuel combustion. This is especially true in a city logistics context \citep{savelsbergh201650th}. 

Over the past years, these observations gave rise to the introduction of pollution- and sustainability-related aspects into traditional Vehicle Routing Problems (VRPs), popularly coined in the literature as the Pollution-Routing Problem \citep{bektacs2011pollution}. In the literature, similar problems and definitions can be found as the Emissions Minimizing VRPs (EMVRPs) \citep{raeesi2019multi} or Green VRPs \citep{erdougan2012green,moghdani2020green}. These models make use of the fact that the amount of transport-related greenhouse gases (GHGs) emissions is directly proportional to the fuel consumption \citep{kirby2000modelling}. Multiple factors are considered, including the slope \citep{suzuki2011new}, vehicle speed \citep{demir2012}, the payload \citep{bektacs2011pollution}, traffic congestion \citep{franceschetti2013time}, driver's operating habit \citep{bandeira2013generating}, and the fleet size and mix \citep{kocc2014fleet}.

Specifically, the PRP aims to build routes that minimize an objective function integrating the vehicle's routing cost (e.g., fuel consumption, the pollution aspect) and the driving costs (e.g., vehicle usage, drivers' wage, and other direct costs aspect).  This paper builds upon this literature and considers modelling the driving costs as fixed costs of vehicles, which is commonly used in heterogeneous vehicle routing problems \citep[see, e.g.,][]{kocc2016thirty}, but often disregarded and handled as duration-dependent costs in PRP. Moreover, considering this rich version of the PRP leads to interesting new methodological challenges for speed optimization.

The majority of models describe the road angle, and thus the network topography, as one of the parameters used to formulate the instantaneous engine-out emission rate \citep{barth2009energy}, but do not consider this further in the model, in the solution methodology or in the results and insights. More specifically, efficient solution methods and extensive computational experiments analyzing the effect of road gradient on fuel consumption and $\text{CO}_2e$ emissions are missing in the literature. One notable exception in the same application domain is the paper by \citet{brunner}. These authors however assume that speed is constant, which is an unrealistic assumption within a city logistics context \citep{franceschetti2013time}. Concerning vehicle speed, several studies are proposing various optimization procedures. However, there is still a need for efficient and fast speed optimization algorithms to use in exact algorithms.

The contributions of this paper are threefold.
\begin{itemize}
    \item We introduce the road gradient to the computation of fuel consumption utilizing terrain elevation information. Despite the inclusion of road gradient into the original formulation of fuel consumption, most papers assume a constant road angle for all vehicle trips.
    \item A number of novel solution approaches are presented, including a branch-and-price algorithm leading to optimal solutions and metaheuristic for larger instances. We develop a novel, fast and efficient algorithm for the vehicle speed optimization. 
    \item We generate important insights based on a broad set of instances to investigate the importance of road gradient on emissions. Additionally, we propose a number of real-life instances.
\end{itemize}

The remainder of this paper is organized as follows. Section \ref{sec: literature review} provides a brief review of related scientific literature. In Section \ref{sec:the-hill-climbing-vehicle-routing-problem}, we give a detailed description of our PRP model. Section \ref{sec: metaheuristics} shows the proposed Tabu Search metaheuristic together with a detailed  description of the fixed-sequence speed optimization algorithm. Later in Section \ref{sec: BCP}, we present the major components of our branch-and-price algorithm. Section \ref{sec:experimental-results} offers extensive computational experiments conducted, and finally, we conclude this research in Section \ref{sec: conclusions}.


\section{Literature Review} \label{sec: literature review}
The significant increasing amount of {$\text{CO}_2e$} emissions derived from road freight transportation has certainly ignited worldwide concerns. Over the last ten years, this resulted in a large body of literature on emissions-aware transportation problems \citep[see, e.g., the survey by][]{moghdani2020green}. The Pollution-Routing Problem (PRP) is an efficient and comprehensive formulation to address the minimization of carbon emissions. The resulting greenhouse emissions of vehicle fuel consumption are the consequence of some influential factors beyond the travel distance \citep{ericsson2001independent,brundell2005influence}. According to \cite{demir2014review}, vehicle fuel consumption is affected by multiple factors such as speed, road gradient, road congestion, driver's operating habit, size and composition (mix) of the vehicle fleet, and payload.

As summarized in Table \ref{PRP-previous-papers}, we observe that some important factors are less studied in previous vehicle routing research, specifically the interaction of load, road gradient, and vehicle speed is missing.

\begin{table}[!h] \centering \small  \caption{Pollution-related Factors Covered by Previous Research on PRPs}\label{PRP-previous-papers}\resizebox{\textwidth}{!}{
\begin{tabular}{lcccclll}
\hline
Authors & \multicolumn{4}{l}{ \makecell[tl]{Pollution-related factors\\ (model formulation)}} & \makecell[tl]{Factors included in the \\ computational analysis} &  \makecell[tl]{Type of \\ dataset} & \makecell[tl]{Solution \\ approach} \\ \cline{2-5}
 & load & speed & slope & others &  &  &  \\ \midrule
\cite{bektacs2011pollution} & $\surd$& $\surd$& $\surd$ &  & load, speed & \begin{tabular}[c]{@{}l@{}}UK\end{tabular} & CPLEX \\
\cite{demir2012} & $\surd$ & $\surd$ & $\surd$ &  & speed & \begin{tabular}[c]{@{}l@{}}UK\end{tabular} & ALNS \\
\cite{franceschetti2013time} & $\surd$ & $\surd$ & $\surd$ & $\surd$ & \begin{tabular}[c]{@{}l@{}}departure time, speed, \\ traffic congestion\end{tabular} & \begin{tabular}[c]{@{}l@{}} UK\end{tabular} & CPLEX \\
\cite{demir2014} & $\surd$ & $\surd$ & $\surd$ &  & driving time, speed & \begin{tabular}[c]{@{}l@{}}UK\end{tabular} & ALNS \\
\cite{kocc2014fleet} & $\surd$ & $\surd$ & $\surd$ & $\surd$ & speed, fleet size and mix & \begin{tabular}[c]{@{}l@{}}UK\end{tabular} & HEA \\
\cite{kramer2015speed} &  & $\surd$ &  & $\surd$ & departure time, speed & \begin{tabular}[c]{@{}l@{}}Modified UK\end{tabular} & ILS \\
\cite{fukasawa2016disjunctive} & $\surd$ & $\surd$ &  &  & speed & \begin{tabular}[c]{@{}l@{}}UK\end{tabular} & DCP \\
\cite{dabia2017exact} &  & $\surd$ & $\surd$ & $\surd$ & start-time, speed & \begin{tabular}[c]{@{}l@{}}UK\end{tabular} & B\&P \\
\cite{rauniyar2019multi} & $\surd$ & $\surd$ &  &  & load & \begin{tabular}[c]{@{}l@{}}UK\end{tabular} & NSGA-II \\
\cite{brunner} & $\surd$ &  & $\surd$ &  & \begin{tabular}[c]{@{}l@{}} arc slope, fixed load \end{tabular} & \begin{tabular}[c]{@{}l@{}}New dataset\end{tabular} & \begin{tabular}[c]{@{}l@{}}Tailored\\ heuristic\end{tabular} \\
\cite{raeesi2019multi} & $\surd$ & $\surd$ & $\surd$ & $\surd$ & \begin{tabular}[c]{@{}l@{}}load, fleet size and mix, \\ road congestion\end{tabular} & \begin{tabular}[c]{@{}l@{}}New dataset\end{tabular} & CPLEX \\
\cite{xiao2020continuous} & $\surd$ & $\surd$ &  & $\surd$ & \begin{tabular}[c]{@{}l@{}}travel-arrival-departure-\\ waiting time, speed, load\end{tabular} & \begin{tabular}[c]{@{}l@{}} UK\end{tabular} & CPLEX \\
\textbf{Proposed work} & $\surd$ & $\surd$ & $\surd$ & $\surd$ & \begin{tabular}[c]{@{}l@{}}\textbf{road gradient}, speed, load\\ \end{tabular} & \begin{tabular}[c]{@{}l@{}}\textbf{Novel}\\ \textbf{dataset}\end{tabular} & \textbf{B\&P, TS} \\ \hline
\end{tabular}}\\
\vspace*{1.5mm}

\textbf{Abbreviations –} CPLEX: IBM Commercial Solver; ALNS: Adaptive Large Neighborhood Search; HEA: Hybrid Evolutionary Algorithm; ILS: Iterated Local Search; DCP: Disjunctive Convex Programming; B\&P: Branch and Price; NSGA-II: Non-dominated Sorting Genetic Algorithm; TS: Tabu Search
\end{table} 

\textbf{Road gradient.} The road grade (slope) has a significant influence on both the conventional-vehicle fuel economy \citep{suzuki2011new} and the electric-vehicle energy consumption \citep{goeke2015routing}. The transit of a typical light-duty vehicle over a sloping road surface, with a +6 percent grade, could increase the fuel consumption by 15-20 percent \citep{boriboonsomsin2009impacts}. The influence of hilly roads is undoubtedly more significant on the fuel consumption of heavy-duty trucks. According to \cite{davis2009transportation}, just a minor increase or decrease in road grade (1-4\%) can reduce or increase fuel economy by more than 50\%. In the case of electric vehicles, the energy consumption is less affected by the road gradient as opposed to conventional vehicles. However, experimental results show the impact of hilly terrain on the electric vehicle miles traveled, confirming that the electric vehicle range in mountainous landscapes is lower \citep{travesset2015transport}.

The road gradient has been mainly studied in PRPs and the so-called Electric Vehicle Routing Problems (E-VRPs). Since the first mathematical formulation of PRPs \citep{bektacs2011pollution}, the road angle was one of the parameters used to define the instantaneous engine-out emission rate. Table \ref{PRP-previous-papers} clearly shows that most PRP-previous contributions included the road angle in the mathematical formulation of fuel use rate. This table also stands out two major gaps in the area that the proposed work is trying to address: (i) despite the inclusion of road gradient into the original formulation of fuel consumption, most of the papers assume that the road angle remains constant throughout all vehicle trips, and (ii) the majority of previously generated problem instances have not yet considered the road gradient.

It is worth mentioning that the present work is the first study in the area to create and optimally solve a set of realistic problem instances, which include elevation information for computing the road slopes along the paths.

As mentioned earlier, we found the paper by \cite{brunner} as the only contribution that studied the influence of road grade on fuel consumption. Although these authors considered the slope in their VRP formulation, they assumed a constant road grade for all arcs that form the directed graph. The above does not reflect accurately the hilly topography profile of a road, which typically connects two nodes (arc) with a sequence of multiple uphill/downhill segments. Another assumption in this contribution is related to the vehicle speed. The authors addressed a VRP without time windows, assuming that the vehicles travel through each arc with a given constant speed (input parameter). Last but not least important with respect to the paper assumptions, in \cite{brunner}, the payload carried by the vehicle is chosen from a prescribed set of values which means that the authors did not consider the payload as a continuous decision variable.

There are very few papers, outside the application context of PRPs, that take into account the effect of road gradients on fuel consumption. For instance, \cite{tavares2009optimisation} utilizing an exponential regression model (COPERT-III method, introduced by \cite{ntziachristos2000copert}) to estimate the minimum fuel consumption during waste collection process. They studied two realistic routing problems showing that the optimal route does not necessarily correspond to the shortest traveled distance. Their computational results demonstrated that significant fuel consumption savings are possible for longer routes with moderate road inclination. Moreover, the contribution of \cite{suzuki2011new} proposed a linear regression model to compute the fuel consumption rate for a heavy-duty truck, based on the fuel-efficiency study of \cite{davis2009transportation}. The author included the road-gradient factor as one of the objective function components (distance and fuel consumption), designed to formulate a traveling salesman problem with time windows.

Concerning the E-VRPs, \cite{yang2014electric} performed a numerical simulation. In their paper, the authors used the electric vehicle’s battery (physical model) theory to study the effects of the road’s slope on electricity consumption for both uphill and downhill paths. They concluded that with the increase of the uphill’s tilt angle, each electric vehicle’s electricity consumption increases significantly. \cite{goeke2015routing} also assumed not-flat terrain with grades in their energy consumption model proposed for electric vehicles. They intensely focused on the effect of load distribution on the performance of commercial electric vehicles. Later, \cite{liu2017impact} also investigate the impact of road gradient on the electricity consumption of electric cars. Using 12 gradient ranges and GPS tracking data with a digital elevation map, the authors showed that, in uphill trips, the energy consumption increases almost linearly with the absolute gradient. However, the numerical results also exhibited the positive effect of regenerative braking power acquired during the downhill trips. Finally, \cite{macrina2019energy} modeled a comprehensive energy consumption function considering the road gradients. Nevertheless, in their computational study, the authors set the road angle equal to zero for conventional and electric vehicles.

\textbf{Vehicle speed.} \cite{bektacs2011pollution} first addressed the speed of the vehicle in the area of PRPs. In the original formulation of this problem, they explicitly assumed that speed over each arc is chosen from a predefined list of possible values. \cite{kocc2014fleet} also utilized the same discrete speed function but analyzed, for the first time, the effect of fleet size and mix in PRPs. The discretization of vehicle speed was also adopted by \cite{eshtehadi2017robust}. Here, the authors considered the demand and travel time uncertainty, both aspects addressed by several robust optimization techniques.
Moreover, \cite{demir2012} proposed a specialized speed optimization algorithm (SOA), which computes optimal speeds on a given path so as to minimize fuel consumption, emissions and driver costs. The authors modified the original SOA, which was initially designed for solving the tramp speed optimization problem \citep{norstad2011tramp}. The speed and traffic congestion were also studied by \cite{franceschetti2013time}, originating the time-dependent PRP. This research considers two phases within the planning horizon, the free-transit phase, and the congested phase. One interesting fact of their computational results is that they reduced the emissions cost by waiting at specific locations (stopped vehicles) and avoiding traffic congestion. \cite{kramer2015speed} also modified the original SOA providing a new speed and departure time optimization algorithm.

Finally, a huge variety of solution methods has been proposed to solve PRPs. One can find the frequent application of metaheuristics such as ALNS \citep{demir2012,demir2014}, evolutionary and genetic algorithms \citep{kocc2014fleet,rauniyar2019multi}, and hybrid approaches \citep{tirkolaee2020multi}. Table \ref{PRP-previous-papers} shows the methods used for solving different variants of PRPs. As can be seen, the utilization of exact algorithms in PRPs is very limited. 
Existing exact methods often approximated or discretized the vehicle speeds in order to reduce the complexity.
\cite{fukasawa2016disjunctive} resolved the issues of speed discretization by introducing a formulation framework to directly incorporate the nonlinear relationship between cost and speed into the PRP. They employed different tools from disjunctive convex programming to find a set of vehicle speeds over the routes, minimizing the total cost (operational and environmental) and respecting the constraints on time and vehicle capacities. More recently, vehicle speed has been computed using continuous optimization on PRP. Here \cite{xiao2020continuous} introduced the continuous PRP ($\epsilon$-CPRP) where the travel time, load flow, departing/arrival/waiting times, and driving speed were treated as continuous decision variables.  The authors developed an $\epsilon$-accurate inner polyhedral approximation method for linearizing the original fuel consumption equation \citep{bektacs2011pollution}, and solved the PRP instances with up to 25 customers.


To the best of our knowledge, \cite{dabia2017exact} is the only previous contribution that addressed a complex variant of PRP, and developed an exact branch-and-price algorithm. To address speed optimization within the pricing subproblem, the authors introduced a ready-time function for updating the speed-dependent routing costs within a bidirectional labeling algorithm. We propose an improved branch-and-price method that uses a novel speed optimization algorithm.

In \cite{dabia2017exact}, a complex ready time function has to be  solved recursively for determining the optimal speed that minimizes the routing cost of a partial path, which can be computationally expensive. In our proposed labeling algorithm, the optimal speed can be determined efficiently and without full ``backtracking'' of the partial path. In addition, we have also developed completion bounds that further improved the efficiency of solving the pricing subproblem.

\section{Problem Formulation}\label{sec:the-hill-climbing-vehicle-routing-problem}

In Section  \ref{sec: co2 emission}, we formulate the carbon dioxide equivalent ($\text{CO}_2e$) emission using the comprehensive modal emissions model (CMEM). In Section 
\ref{sec: model}, we introduce a mixed-integer linear programming formulation for the joint Pollution Routing and Speed Optimization problem (PRP-SO).

\subsection{Modelling $\text{CO}_2e$ emissions} \label{sec: co2 emission}
We model $\text{CO}_2e$ emissions using the CMEM (refers to e.g. \citet{barth2009energy,boriboonsomsin2009impacts, demir2012}). The parameters and the values for light, medium, and heavy-duty vehicles (denoted respectively as LDV, MDV, HDV) we used in our experiments are shown in Table \ref{CMEM parameters}. 
\begin{table}[!ht]\small  \caption{Comprehensive Modal Emissions Model (CMEM)}\label{CMEM parameters}
\centering
\begin{tabular}{clccc}
\toprule
Symbol & Description & LDV  & MDV  & HDV  \\ 
\midrule
$F_k$ & Engine friction factor (kJ/rev/liter) & 0.23 & 0.20 & 0.17 \\ 
$N_k$ & Engine speed (rev/s) & 35 & 34 & 33 \\ 
$V_k$ &  Engine displacement (liters) & 3 & 7 & 11 \\ 
$A_k$ &  Frontal surface area of a vehicle (m$^2$) & 5  & 7.6 & 8.2\\ 
$C^d_k$ & Aerodynamic drag coefficients  & 0.32 & 0.55 & 0.70\\ 
$C^r_k$ & Rolling resistance coefficients & 0.01 & 0.009 & 0.008\\ 
$r_k$ & Vehicle acceleration (m$/$s$^2$) & 0 & 0 & 0 \\ 
$w_k$ & Curb weight (kg) &2,300 &5,500 &13,000 \\
$\kappa$  & Heating value for diesel fuel (kJ/g) & 45 & 45 & 45\\ 
$\varepsilon$ & Vehicle drive train efficiency  & 0.4 & 0.4 & 0.4 \\ 
$\varpi$ & Efficiency parameter for diesel engines  & 0.9 & 0.9 & 0.9\\ 
$\xi$ & Fuel-to-air mass ratio & 1 & 1 & 1\\ 
$\psi$ & Conversion factor from grams to liters & 737 & 737 & 737 \\ 
$\rho$ &  Air density (kg/m$^3$) & 1.2041 & 1.2041 & 1.2041 \\ 
$g$ & Gravity (m$/$s$^2$)  & 9.81  & 9.81 & 9.81\\
\bottomrule
\end{tabular} 
\end{table}

According to the CMEM model, the instantaneous fuel use rate of a vehicle $k$ when traveling at a constant speed $\nu_k$ with payload $x$ on a path with the road angle $\phi$ is given by
$$\frac{\xi}{\kappa\psi }\Big(  F_k N_k V_k + \frac{0.5C^d_k A_k\rho\nu_k^3 + (w_k +x)\nu_k (r_k + g sin\ \phi + gC^r_k\ cos\ \phi) }{1000  \varepsilon  \varpi }\Big) $$
When traversing a distance of $d$ meters, the amount of fuel consumption is therefore given by
$$\frac{\xi F_k N_k V_k}{\kappa \psi } \frac{d}{\nu_k} + \frac{1}{1000  \varepsilon  \varpi }(r_k + g\ sin\ \phi + g\ C^r_k\ cos\ \phi)(w_k+x)(d) + \frac{0.5C^d_k A_k \rho }{1000  \varepsilon  \varpi } d\nu_k^2 $$

Now that, there are a sequence of road segments associated on each arc $e$ which are denoted by $S_e$. Let $d_{es}$ denote the travel distance of segment $s\in S_e$, and $x_{ke}$ denote the payload of vehicle $k$ when traversing arc $e$. For traversing the sequence of road segments associated on arc $e$, the amount of fuel consumption of vehicle $k$ can then be determined by
\begin{align}
\alpha_{ke}\frac{1}{\nu_k} + \beta_{ke}(w_k + x_{ke}) + \gamma_{ke} \nu_k^2  \label{co2 cost} 
\end{align}

where 
\begin{align}
&\alpha_{ke} =\frac{\xi F_k N_k V_k \sum_{s \in S_e}d_{es}}{\kappa \psi } \\
&\beta_{ke} =\frac{\xi \sum_{s \in S_e} d_{es} (r_k + g\ sin\ \phi_{es} + g\ C^r_k\ cos\ \phi_{es})}{1000  \varepsilon  \varpi \kappa\psi} \label{co2 beta} \\
&\gamma_{ke}=  \frac{0.5 \xi C^d_k A_k \rho \sum_{s \in S_e} d_{es}}{1000  \varepsilon  \varpi \kappa \psi}
\end{align}
To speed up the $\text{CO}_2e$ emission calculations, we precompute the parameters $\alpha_{ke}$, $\beta_{ke}$, and $\gamma_{ke}$ for all the vehicles $k$ and arcs $e$.
These parameters are used in formulating the mathematical model in Section \ref{sec: model} and developing the solution approaches in Section \ref{sec: metaheuristics} and Section \ref{sec: BCP}.


\subsection{Mathematical model} \label{sec: model}

\begin{table}[!ht]\small  \caption{Notation for the MINLP Model}\label{MIP notation} \centering
\begin{tabular}{cll}
\toprule
Symbol & Description & Domain \\ 
\midrule
$\cal K$ & Set of vehicles & index set \\ 
$\cal N$  & Set of nodes representing $n$ customers and the depot& index set \\ 
$\cal A$  & Set of arcs & index set \\ 
$q_i$ & Demand of customer $i$  & $\mathbb{Z}^+$  \\ 
$w_k$ & Curb weight of vehicle $k$& $\mathbb{Z}^+$  \\ 
$Q_k$ & Capacity of vehicle $k$  & $\mathbb{Z}^+$  \\ 
$e_i$ & Earliest start time at customer $i$  & $\mathbb{R}^+$  \\ 
$l_i$ & Latest start time at customer $i$  & $\mathbb{R}^+$  \\  
$s_i$ & Service time at  customer $i$ & $\mathbb{R}^+$  \\  
$f_k$ & Fixed cost of vehicle $k$& $\mathbb{R}^+$  \\ 
$c_k$ & Variable cost of vehicle $k$ including fuel and $\text{CO}_2e$ emission costs & $\mathbb{R}^+$ \\ 
$a_k$ & Lower limit of the speed of vehicle $k$  & $\mathbb{R}^+$  \\ 
$b_k$ & Upper limit of the speed of vehicle $k$  & $\mathbb{R}^+$  \\ 
$d_e$ & Distance of arc $e$ & $\mathbb{R}^+$  \\
$\alpha_{ke}$ & A constant for estimating the fuel consumption of vehicle $k$ on arc $e$ & $\mathbb{R}^+$ \\
$\beta_{ke}$ & A constant for estimating the fuel consumption of vehicle $k$ on arc $e$ & $\mathbb{R}^+$ \\
$\gamma_{ke}$ & A constant for estimating the fuel consumption of vehicle $k$ on arc $e$ & $\mathbb{R}^+$ \\
$z_{ki}$ & Decision variable, to allocate customers to vehicles & $\mathbb{B}$\\ 
$y_{ke}$ & Decision variable, to allocate arcs to vehicles & $\mathbb{B}$\\ 
$x_{ke}$ & Decision variable, payload of vehicle $k$ when traversing on arc $e$ & $\mathbb{R}^+$\\ 
$t_{ki}$ & Decision variable, time at which vehicle $k$ starts serving customer $i$ & $\mathbb{R}^+$\\ 
$\nu_{k}$ & Decision variable, speed of vehicle $k$ & $\mathbb{R}^+$\\ 
\bottomrule
\end{tabular} 
\end{table}

Let $\cal K$ be the set of vehicles. Let $G({\cal N}, {\cal A})$ be the underlying directed graph. A set of nodes  ${\cal N} =\{0, 1, ...,n\}$ contains $n$ customers and a depot (represented by node 0). The set of arcs  ${\cal A}$, defined as $\{ (i,j) \in {\cal N} \times {\cal N}: i \not= j \}$, represents the paths between the nodes. Each node $i \in \cal N$ is associated with a demand $q_i$ and a time window [$e_i$, $l_i$]. Each vehicle $k \in  \cal K$ is associated with a fixed cost $f_k$, a variable cost $c_k$ (including the costs for $\text{CO}_2e$ emission and fuel consumption), a vehicle capacity $Q_k$, vehicle curb weight $w_k$, and speed limits [$a_k$, $b_k$].
Each arc $e \in {\cal A}$ is associated with a distance $d_e$, and the parameters $\alpha_{ke}$, $\beta_{ke}$ and $\gamma_{ke}$ described in Section \ref{sec: co2 emission} for estimating fuel consumption.

For all $k \in \cal K$ and $i \in {\cal N}$, let $z_{ki}$ be a binary decision variable, with $z_{ki}=1$ if and only if customer $i \in {\cal N} \setminus \{0\}$ is served by vehicle $k$; and with $z_{k0}=1$ if and only if vehicle $k$ is in use. For all $k \in \cal K$ and $e \in {\cal A}$, let $y_{ke}$ be a binary decision variable with $y_{ke}=1$ if and only if vehicle $k$ traverses arc $e$, and let $x_{ke}$ denote the corresponding payload when vehicle $k$ traverses arc $e$.
For all $k \in \cal K$ and $i \in {\cal N}$, let $t_{ki}$ denote the time at which vehicle $k$ starts serving customer $i \in {\cal N} \setminus \{0\}$; and let $t_{k0}$ denote the time vehicle $k$ returns to the depot.
For all $k \in \cal K$, let $\nu_{k}$ denote the speed of vehicle $k$.

The objective is to minimize the total $\text{CO}_2e$ emission costs, fuel costs, and vehicle fixed costs, subject to the following constraints: 
i) the total demand in a vehicle does not exceed the vehicle capacity; 
ii) every route starts and ends at the vehicle's home depot;
iii) every customer is visited exactly once by exactly one vehicle;
iv) all vehicles should return to its home depot within a time limit;
v) every vehicle travels at a speed within the speed limits.

The joint Pollution Routing and Speed Optimization problem (PRP-SO) can be formulated as the following nonlinear integer programming model.
\begin{align}
(\text{PRP-SO}):\nonumber  \\
\ \text{min}\ &\sum\limits_{k \in {\cal K}} f_k z_{k0} + \sum\limits_{k \in {\cal K}} \sum\limits_{e \in {\cal A}}  c_k y_{ke}\Big(\frac{\alpha_{ke}}{\nu_k} + \beta_{ke}w_k + \beta_{ke}x_{ke} + \gamma_{ke} \nu_k^2 \Big), \label{MIP: obj}  \\
\text{s.t.}\ \ 
&\sum\limits_{k \in {\cal K}} z_{ki} = 1, &&\forall i \in {\cal N} \setminus \{0\}, \label{MIP: visit once}\\
&\sum\limits_{e \in \delta^+(i)} y_{ke} = \sum\limits_{e \in \delta^-(i)} y_{ke} = z_{ki}, &&\forall k \in {\cal K}, i \in {\cal N}, \label{MIP: y-flow conservation} \\
&\sum\limits_{e \in \delta^-(i)} x_{ke} - \sum\limits_{e \in \delta^+(i)} x_{ke} = q_i z_{ki}, &&\forall k \in {\cal K}, i \in {\cal N} \setminus \{ 0\}, \label{MIP: x-flow conservation} \\
&x_{ke} \leq (Q_k - q_i) y_{ke},&&\forall k \in {\cal K}, e=(i,j) \in \cal A, \label{MIP: arc capacity}\\
&t_{kj} - t_{ki} \geq s_i + d_e\ \frac{y_{ke}}{\nu_k}   - l_0 \Big(1- y_{ke}\Big), &&\forall  k \in {\cal K}, e=(i,j) \in {\cal A}: j\neq 0,\label{MIP: time windows 1} \\
&  e_{i}z_{ki} \leq  t_{ki} \leq  l_{i}z_{ki}, &&\forall  k \in {\cal K}, i \in {\cal N}\setminus \{0\},\label{MIP: time windows 1.5}  \\
& t_{ki} + s_i + d_e\ \frac{y_{ke}}{\nu_k}  \leq l_0, &&\forall  k \in {\cal K}, i \in {\cal N}\setminus \{0\}, e=(i,0)\label{MIP: time windows 2}  \\
&a_k \leq \nu_k \leq b_k,&&\forall k \in {\cal K}, \label{MIP: speed}  \\
&\nu_k \in \mathbb{R}^+, &&\forall k \in {\cal K},\label{MIP: speed variables} \\
&t_{ki} \in \mathbb{R}^+, &&\forall  k \in {\cal K}, i \in  {\cal N}, \\
&x_{ke} \in \mathbb{R}^+, &&\forall k \in {\cal K}, e \in {\cal A},\label{MIP: arc flow variables}  \\
&y_{ke} \in \{0, 1\}, &&\forall k \in {\cal K}, e \in {\cal A},\label{MIP: arc binary variables}  \\
&z_{ki} \in \{0, 1\}, &&\forall k \in {\cal K}, i \in  {\cal N}.\label{MIP: arc variables}  
\end{align}

The objective function \eqref{MIP: obj} minimizes the total vehicle fixed costs, fuel consumption costs, and $\text{CO}_2e$ emission costs.
Constraints \eqref{MIP: visit once} - \eqref{MIP: y-flow conservation}  ensure that each customer is visited once by exactly one vehicle.
Constraints \eqref{MIP: x-flow conservation} are the flow conservation constraints.
Constraints \eqref{MIP: arc capacity} ensure that the payload does not exceed vehicle capacity.
Constraints \eqref{MIP: time windows 1} - \eqref{MIP: time windows 2} ensure that customers are visited within the given time windows.
Constraints \eqref{MIP: speed} ensure that vehicles travel at a speed within the limits.

\section{Metaheuristic} \label{sec: metaheuristics}
Metaheuristic approaches require frequently evaluating solutions with fixed vehicle routes. We present an efficient algorithm for finding the optimal vehicle speed of a given  route so that we can compute the costs of $\text{CO}_2e$ emissions and fuel consumption efficiently.
In Section \ref{sec: Fixed-sequence speed optimization},  we present a novel polynomial-time algorithm for solving the fixed-sequence speed optimization subproblem --- determine the optimal speed of a vehicle when the customer sequence is fixed. To demonstrate the effectiveness of the speed optimization algorithm, it is embedded into a Tabu Search (TS) metaheuristic for our experiments. TS is originally proposed by \cite{glover1986future} as a synthesis of the perspectives of operations research and artificial intelligence. TS has been widely used and has been shown to be effective for finding near-optimal solutions to vehicle routing problems, see, e.g. \cite{gendreau1994tabu}, \cite{toth2003granular}, \cite{lai2016tabu}. Review on the recent development of TS can refer to \cite{GendreauTabuSearch}.
Major components of the proposed TS include the fixed-sequence speed optimization subproblem in Section \ref{sec: Fixed-sequence speed optimization}, the penalized objective function in Section \ref{sec: objective function}, the initial solutions in Section \ref{sec: initial solution}, the neighborhood structure in Section \ref{sec: neighborhood structure}, the intra-route improvement procedure in Section \ref{subsec:intra-route-improvement-procedure}, and the search procedure in Section \ref{sec: procedure}.

\subsection{Fixed-sequence speed optimization subproblem}  \label{sec: Fixed-sequence speed optimization}
The fixed-sequence speed optimization subproblem determines the optimal speed of a vehicle for a given customer sequence.  Let $R=(v_1, v_2, ..., v_n)$ denote a fixed sequence of nodes in a vehicle route where both $v_1$ and $v_{n}$ represent the home depot.

Without time-window constraints, vehicles should always travel at a speed that is most cost-efficient and within the speed limits of the vehicle.
The most cost-efficient speed of vehicle $k$ when there are no time window constraints is given by
$$\bar v_k= \sqrt[3]{\frac{\xi_k F_k N_k V_k 1000  \varepsilon_k  \varpi_k }{ C^d_k A_k \rho_k }}$$
Since $\bar v_k$ is independent of the vehicle route when there are no time window constraints, the optimal speed of vehicle $k$ within speed limits $[a_k, b_k]$ can be computed by
$$v^*_k = \begin{cases}
              \bar v_k,  &\text{if $a_k \leq \bar v_k \leq b_k$},\\
              a_k, &\text{if $\bar v_k \leq a_k$},\\
              b_k, &\text{if $\bar v_k \geq b_k$}.
\end{cases}$$

In the presence of time window constraints, a vehicle should travel at a speed that satisfies all time window constraints and incurs a minimal cost.
Let $\sigma(R)$ denote the lowest vehicle speed without violating any time windows associated on the nodes of route $R$.
When a vehicle travels at the speed of $\sigma(R)$ in route $R$, all time window constraints will be satisfied.
We will show in Proposition \ref{sec: Algorithm for determining the lowest feasible speed} that $\sigma(R)$ can be computed efficiently by
\begin{align}
    \sigma(R) = \max_{i,j \in \{1,2,...,|R|\}:i< j} \frac{\Delta(j)-\Delta(i)}{l_{v_j}-e_{v_i}-S_{ij}} \label{fixed-seq optimal speed 1}
\end{align}
where $\Delta(i)=\sum_{l=1}^{i-1} d_{v_l, v_{l+1}}$ denotes the total distance from the depot to node $i$ along vehicle route $R$, $S_{ij} = \sum_{k=i}^{j-1}  s_{v_k}$ denotes the total service time from node $v_{i}$ to node $v_{j-1}$ along vehicle route $R$, and that $[e_v, l_v]$ is the time window associated on node $v$.
If $\sigma(R) \leq 0$, no feasible speed exists due to conflicting time window constraints.

Since a vehicle can wait at the customer node if the vehicle arrives earlier than the lower bound of a time window, any speed that is greater than $\sigma(R)$ also satisfies the time-window constraints in route $R$. Thus, the optimal speed of vehicle $k$ when traversing on route $R$ is given by
\begin{align}
    \max(v^*_k,\ \sigma(R)) \label{fixed-seq optimal speed 2}
\end{align}
The total cost (as defined in  the objective function \eqref{MIP: obj}) of a given solution can be evaluated straightforward when the optimal speeds of the vehicle routes have been found by using \eqref{fixed-seq optimal speed 1} and \eqref{fixed-seq optimal speed 2}.

\begin{proposition} \label{sec: Algorithm for determining the lowest feasible speed}
The minimum speed without violating any time window constraints of a vehicle route $R$ is given by $\sigma(R)$ when $\sigma(R)>0$, and can be obtained in a $O(n^2)$ complexity where $n$ is the number of nodes in route $R$.
\end{proposition}
\proof{} 
    Let $R=(v_1, v_2, ..., v_n)$ denote the sequence of nodes in a vehicle route $R$ with both $v_1$ and $v_{n}$ representing the home depot, and there is at least one customer node in $R$ i.e. $n\geq 3$.
Furthermore, we assume that the nodes in $R$ are not all located at the same location.
  For $i=2,...,n$, let $\Delta(i)=\sum_{l=1}^{i-1} d_{v_l, v_{l+1}}$ denote the total distance from node $v_1$ (the depot) to node $v_i$ along vehicle route $R$. Set $\Delta(1)=0$.
    For all $i=1,2,...,n$ and $j=2,...,n$, with $i<j$, let $S_{ij} = \sum_{k=i}^{j-1}  s_{v_k}$ denote the total service time from node $v_{i}$ to node $v_{j-1}$.
    For every distinct pair of nodes $i,j \in \{1,2,...,n\}$, with $i< j$, the following lower bound of $\sigma(R)$ can be derived:
    $$\frac{\Delta(j)-\Delta(i)}{l_{v_j}-e_{v_i}-S_{ij}}$$
    Any vehicle speed that is higher than the above lower bound induced from nodes $i$ and $j$ would violate one of the time window constraints associated on nodes $v_i$ and $v_j$.

    The minimum speed without violating any time window constraints in $R$ is given by the maximum lower bounds for all distinct pairs of nodes, and thus we have
    $$\sigma(R) = \max_{i,j \in \{1,2,...,|R|\}:i< j} \frac{\Delta(j)-\Delta(i)}{l_{v_j}-e_{v_i}-S_{ij}}.$$
    Since there are $\frac{n (n - 1)}{2}$  distinct pairs of nodes in a vehicle route of length $n$, $\sigma(R)$ can be computed in $O(n^2)$.

    To complete the proof, we will show that a pair of conflicting time window constraints exist in route $R$ if and only if $\sigma(R) \leq 0$.
    Consider two cases: i) $\sigma(R) < 0$; ii) $\sigma(R) = 0$.
    Case i: By definition we have $\Delta(j)-\Delta(i)\geq 0$ for all $i,j=1,...,n$ with $i< j$.
    Therefore, $\sigma(R) < 0$ iff there exist $v_i$ and $v_j$ with $l_{v_j}-e_{v_i}< S_{ij}$ which implies a violation on one of the time window constraints associated on nodes $v_i$ and $v_j$.
    Case ii: Suppose the contrary, we consider a vehicle with a speed equals to zero and that all the time window constraints in $R$ are satisfied.
    Since the due date $l_{v} \neq \infty$ for all $v$, we have $\sigma(R) = 0$ iff $\Delta(j)-\Delta(i)=0$ for all $i,j=1,...,n$ with $i< j$.
    Since the nodes in $R$ are not all located at the same location, there exist two distinct nodes $v_i$ and $v_j$ in $R$ with $d_{v_i, v_{j}}>0$ and $i< j$.
    This implies that, there exist $v_i$ and $v_j$ in $R$ with $i< j$ and $\Delta(j)-\Delta(i)>0$ which leads to a contradiction.
\Halmos
\endproof

\subsection{Penalized objective function} \label{sec: objective function}
We allow infeasible solutions in the search space.
It is implemented as a penalized objective function that is obtained by relaxing some of the constraints and incorporating them into the objective function with the use of self-adjusting penalty parameters.

If a vehicle $k \in \cal K$ is assigned to a non-empty route $R$, and is travelling at a speed of $\nu_k$, the travel cost $c(R)$ can be written as
$c(R)=f_k + \sum_{e \in {\cal A}(R)}  c_k \big(\frac{\alpha_{ke}}{\nu_k} + \beta_{ke}w_k + \beta_{ke}\bar x_{ke} + \gamma_{ke} \nu_k^2 \big) $
where $\bar x_{ke}$ is the payload of vehicle $k$ on arc $e$ and ${\cal A}(R)$ denotes the arcs on route $R$.
The overload $\text{P}(R)$ representing the violation of the vehicle capacity constraints is defined as
$\text{P}(R) = \big[  \sum_{ i \in  {\cal N}(R)} q_i - Q_k  \big]^+$ where ${\cal N}(R)$ denotes the customer nodes on route $R$.
After incorporating the penalties for possible violations, the penalized objective function value is computed by $z(R) = c(R) + \rho \text{P}(R) $ where $\rho \in \mathbb{R}^+$ is the penalty weight that is self-adjusting in the search.
The penalty weight $\rho$ is initialized as 1, and updated in every $\delta$ iterations as follows.
If $\text{P}(R)=0$ for all vehicle routes $R$, then update $\rho$ to $0.5\rho$; otherwise, update $\rho$ to $2\rho$.


\subsection{Initial solution} \label{sec: initial solution}
The initial solution is created by first applying a stochastic insertion-based heuristic and then improved by using the TS procedure described in Algorithm \ref{algo: Tabu search} with a limited number of iterations.
The number of iterations is set to $\lceil I_1 n \rceil$ where $I_1$ is a user-controlled parameter and $n$ is the number of customers.
After creating ten such initial solutions, the best one is returned as the initial solution for the main search procedure described in Section \ref{sec: procedure}.

The following stochastic insertion-based heuristic is applied.
It begins by assigning exactly one arbitrarily selected customer to each of the randomly selected $l$ vehicles, where $l$ is a randomly generated integer between 1 and the total number of vehicles available.
Then, the remaining customers are considered one by one, following a randomized order.
For each customer, all possible locations in all routes are evaluated for insertion, and the customer is subsequently inserted into a vehicle route at the position that minimizes the insertion cost (the incremental change in the penalized objective function value).
We determine the insertion costs by using the algorithm described in Section \ref{sec: Fixed-sequence speed optimization} that performs speed optimization for a fixed sequence of nodes on a route.

\subsection{Neighborhood structure} \label{sec: neighborhood structure}
Let $\mathcal X$ denote the set of all feasible solutions for the instance.
We define the solutions in the neighborhood of a given solution $ x \in \mathcal X$ as $\textbf{N}(x)$.
At each iteration, all possible move operations for all customers are evaluated and the best one is subsequently performed.
The move operation involves relocating a customer from its current route to another route at the location that minimizes the insertion cost.
The insertion cost is determined by applying the speed optimization algorithm described in Section \ref{sec: Fixed-sequence speed optimization}.
The best move operation is the one that leads to the lowest total penalized objective function value and the following diversification penalty.
For a given solution $x \in \mathcal X$, we define the diversification penalty as $\phi(x) =  \lambda c(x) \sqrt {n} \vartheta_{ir} $ where $n$ is the number of customers, $\vartheta_{ir}$ counts the number of times customer $i$ has been moved to route $r$ so far in the search, and $\lambda$ is a positive parameter that controls the intensity of diversification.
Readers can refer to  \cite{soriano1996diversification} for extensive discussion on diversification schemes.

To prevent cycling, if a customer has been moved from route $r$ to route $s$ in a given iteration, then moving the same customer back to route $r$ is declared tabu which implies that this reverse operation is forbidden for the next $\lceil h\ \text{log}_{10}(n) \rceil$ iterations where $h$ is a user-controlled parameter and $n$ is the number of customers.
To prevent the search from stagnating, the following aspiration criterion is applied: a Tabu move is allowed only when the resulting solution is feasible and has an objective function value that is better than that of the current best feasible solution found by the search.

\subsection{Intra-route improvement procedure}\label{subsec:intra-route-improvement-procedure}
The following procedure is applied for improving the solution by modifying customers' position within the same route: a customer is randomly picked and then reinserted into the best location of the same route, until no further improvement is possible. The intra-route improvement procedure is invoked after the selected move operation is performed, the parameter for overload penalty is updated, or when the best feasible solution is improved.

\subsection{Search procedure} \label{sec: procedure}
The TS routine is presented in Algorithm \ref{algo: Tabu search}. The search procedure consists of two phases.
The first phase of the procedure starts by constructing an initial solution as described in Section \ref{sec: initial solution}. Next, the best solution found in the first phase is improved in the second phase by executing $I_2$ iterations of the TS routine.

\begin{algorithm}[!ht]
    \caption{Tabu Search}\label{algo: Tabu search}
    \begin{algorithmic}[1]
        \State \textbf{input:} initial solution $x_0$
        \State Set $x=x_0$. If $x$ is feasible, set $z^*=c(x)$ and $x^*=x$; otherwise, set $z^*=\infty$ and $x^*=x$.
        \State determine $z (\bar x)$ and $\phi(\bar x)$ for all $\bar x \in {\cal N}(x)$
        \While{stopping condition is not satisfied}
            \State  \textbf{select} $\bar x \in {\cal N}(x)$ \textbf{that}
            \State - minimizes $z(\bar x)$ + $\phi(\bar x)$
            \State - $\bar x$ is non-tabu or it satisfies the aspiration criteria
            \State set the reverse move tabu for $\theta$ iterations
            \State perform the intra-route improvement procedure on $\bar{x}$. \If{$\bar x$ is feasible and $z(\bar{x}) < z^*$} set $x^*={\bar x}$ and $z^*=z(\bar{x})$.
            \EndIf
            \State set $x=\bar{x}$, and update the penalty weight of overload for every $\delta$ iterations
            \State update $z (\bar x)$ and $\phi(\bar x)$ for all $\bar x \in \textbf{N}(x)$
        \EndWhile
        \State \textbf{return} $x^*$
    \end{algorithmic}
\end{algorithm}

Algorithm \ref{algo: Tabu search} shows the TS procedure.
In line 2, the algorithm starts with an initial solution $x_0$ which is obtained by running the TS routine with $I_1$ iterations as described in Section \ref{sec: initial solution}. In line 3, the penalized objective function and diversification penalty associated on all the move operations are updated. The loop in lines 4--12 is then invoked and stops until after $I_2$ iterations. In line 5--7, the best neighborhood solution is picked which takes into account the diversification penalty, overload penalty, solutions in the Tabu list, and the aspiration criteria.
As described in Section \ref{sec: neighborhood structure}, the diversification penalty $\phi(x)$ depends on the intensification parameter $\lambda$. In line 8, the reverse move is forbidden for the next $\lceil h\ \text{log}_{10}(n) \rceil$ iterations. As shown in lines 9--10, the intra-route improvement procedure described in Section \ref{subsec:intra-route-improvement-procedure} is performed on the two routes that are modified in the best neighborhood solution. The incumbent is updated if a new best feasible solution is identified.
In lines 11--12, the search moves to the selected neighboring solution. The penalty weight of overload is updated in every $\delta$ iterations. Whenever the search moves to a neighboring solution, the penalized objective function and diversification penalty associated on each of the move operations are again updated by using the speed optimization algorithm described in Section
\ref{sec: Fixed-sequence speed optimization} and the penalized objective function defined in Section \ref{sec: objective function}. The computational time is manageable since we only need to update the costs of the move operations that are involved with the modified routes.
The values for the algorithmic parameters in the search procedure include $I_1$, $I_2$, $\lambda$, $h$ and $\delta$ which are set according to the parameter tuning experiment described in Section \ref{subsec:-parameter-tuning}.


\section{Branch-and-price Algorithm} \label{sec: BCP}
Branch-and-price is a successful exact method for solving several VRP variants. It relies on reformulating the problem as a set-partitioning (SP) problem and applying branch-and-bound to solve the SP formulation, which is known to provide tight linear bounds. As the number of variables in the SP model grows exponentially with the instance size, column generation is applied to identify profitable variables and add them to the model dynamically. For a review of branch-and-price applied to vehicle routing we refer to \cite{costa2019exact}. In this section, we present the SP-based PRP-SO formulation and describe our branch-and-price solution method, with a focus on the specialized algorithm for solving the non-linear column generation subproblem.

\subsection{Set-partitioning formulation}\label{sec:spf}
A route $p$ is feasible for vehicle $k$ when i) there exists a value $\nu\in[a_k,b_k]$ such that no time-window constraint is violated when vehicle $k$ executes route $p$ at speed $\nu$; and ii) the total demand of the customers served along $p$ does not exceed the vehicle capacity $Q_k$. Let $\Omega_k$ be the set of feasible routes for vehicle $k$. Furthermore, for each $p\in\Omega_k$, let $c_p^k$ be the cost incurred when vehicle $k$ executes route $p$ at the optimal speed, as defined in Section \ref{sec: Fixed-sequence speed optimization}. Finally, for each route $p\in\Omega_k$ and customer $i\in{\cal N}\setminus\{0\}$, let 
$$\delta_{ip}=\begin{cases}
1, &\text{if route $p$ visits customer $i$,}\\
0, &\text{otherwise.}
\end{cases}$$

We define binary decision variables $\lambda_{p}^k$ for all $k\in\cal K$ and $p\in\Omega_k$, such that $\lambda_{p}^k=1$ if and only if vehicle $k$ executes route $p$ in the solution. Then, the PRP-SO is formulated as the following SP problem:
\begin{align}
(\text{SP}):\ \min\ &\sum\limits_{k \in {\cal K}}\sum\limits_{p \in \Omega_k} c_p^k \lambda^k_p, \label{sp: obj} \\
\text{s.t.} \ & \sum\limits_{k \in {\cal K}}\sum\limits_{p \in \Omega_k} \delta_{ip}\lambda_p^k = 1, &&\forall i \in {\cal N} \setminus \{0\}, \label{sp: visit} \\
&\sum\limits_{p \in \Omega_k}  \lambda^k_p \leq 1, &&\forall k \in \cal K,  \label{sp: fleetsize} \\
&\lambda_p^k \in \{0,1\}, &&\forall k \in {\cal K}, \forall p \in \Omega_k.\label{sp: domain} 
\end{align}

The objective function \eqref{sp: obj} minimizes the total costs. Constraints \eqref{sp: visit} ensure that each customer is visited exactly once, and constraints \eqref{sp: fleetsize} enforce a maximum of one route per vehicle.

\subsection{Column generation subproblem}\label{sec:colgen}
The restricted master problem (RMP) refers to the linear relaxation of \eqref{sp: obj}-\eqref{sp: domain} with a restricted set of vehicle routes. Let $\pi_i$, $i\in{\cal N}\setminus \{0\}$, be the dual prices associated with constraints \eqref{sp: visit} after solving the RMP. The pricing problem for vehicle $k$ is defined as follows:
\begin{align}
\text{($\text{PP}$)}:\ \min\ &f_k+c_k \sum\limits_{e\in{\cal A}}\Big(\frac{\alpha_{ke}}{\nu}+\beta_{ke}(w_k+x_{e})+\gamma_{ke}\nu^2\Big)y_e-\sum\limits_{i\in{\cal N}\setminus\{0\}}\pi_i z_{i}, \label{sub: obj}\\
\text{s.t}.\ \ 
&\sum\limits_{e \in \delta^+(i)} y_{e} = \sum\limits_{e \in \delta^-(i)} y_{e} = z_{i}, &&\forall i \in {\cal N} \setminus \{0\}, \label{sub: y-flow conservation 1} \\
&\sum\limits_{e \in \delta^+(0)} y_{e} = \sum\limits_{e \in \delta^-(0)} y_{e} = 1, \label{sub: y-flow conservation 2} \\
&\sum\limits_{e \in \delta^-(i)} x_{e} - \sum\limits_{e \in \delta^+(i)} x_{e} = q_i z_{i}, &&\forall i \in {\cal N} \setminus \{ 0\}, \label{sub: x-flow conservation} \\
&x_{e} \leq (Q_k - q_i) y_{e},&&\forall e=(i,j) \in \cal A, \label{sub: arc capacity}\\
&\sum\limits_{i \in {\cal N}\setminus \{0\}}q_i z_{i} \leq Q_k,  \label{sub: capacity} \\
&t_{j} - t_{i} \geq s_i + d_e\ \frac{y_{e}}{\nu}   - l_0 \Big(1- y_{e}\Big), &&\forall e=(i,j) \in {\cal A}: j\neq 0,\label{sub: time windows 1} \\
&e_{i}z_{i} \leq  t_{i} \leq l_i z_{i}, &&\forall  i \in {\cal N}\setminus \{0\},\label{sub: time windows 1.5}  \\
& t_{i} + s_i + d_e\ \frac{y_{e}}{\nu}  \leq l_0, &&\forall i \in {\cal N}\setminus \{0\}, e=(i,0),\label{sub: time windows 2}  \\
&\nu\in[a_k,b_k],&& \label{sub: speed}  \\
&t_{i} \in \mathbb{R}^+, &&\forall  i \in  {\cal N}, \\
&x_{e} \in \mathbb{R}^+, &&\forall e \in {\cal A},\label{sub: arc flow variables}  \\
&y_{e} \in \{0, 1\}, &&\forall e \in {\cal A},\label{sub: arc binary variables}  \\
&z_{i} \in \{0, 1\}, &&\forall i \in  {\cal N} \setminus\{0\}.\label{sub: arc variables}  
\end{align}

The objective function \eqref{sub: obj} minimizes the route-dependent costs (including vehicle fixed cost and $\text{CO}_2e$ emissions cost) and the dual prices of the RMP. Constraints \eqref{sub: y-flow conservation 1}-\eqref{sub: y-flow conservation 2} are the flow conservation constraints. 
Constraints \eqref{sub: x-flow conservation}-\eqref{sub: arc capacity} keep track of the payload in the vehicle along each arc traversed.
Constraint \eqref{sub: capacity} is the vehicle capacity constraint. Constraints \eqref{sub: time windows 1}-\eqref{sub: time windows 2} are the time window constraints. Finally, constraint \eqref{sub: speed} ensures that the vehicle travels at a speed within the prescribed limits.

\subsection{Pricing algorithm} \label{sec: pricing}
The pricing algorithm identifies columns with a negative reduced cost, that is, it finds solutions to (PP) with a negative objective value. As common in branch-and-price for vehicle routing, we solve the pricing problem with a labeling algorithm. A label represents a partial route from the depot to a customer. Label extensions are created by extending labels to all feasible customers. Table \ref{pricing:labels} summarizes the attributes of a label alongside their corresponding initialization values and updating rules.

\begin{table}[!h]\small  \caption{Pricing Algorithm: Label Attributes, Initialization Values and Updating Rules}\label{pricing:labels}
\begin{minipage}{\textwidth}
\centering
\begin{tabular}{llll}
\toprule
Attribute$^a$ & Description & Initialization$^b$ & Updating rule$^c$ \\
\midrule
$N_P$ & Set of customers & $\{i\}$ & $N_Q=N_P\cup\{j\}$\\
$n_P$ & Last customer served & $i$ &$j$ \\
$q_P$ & Total demand & $q_i$ & $q_Q=q_P+q_j$\\
$D_P$ & Total distance traveled & $d_e$ & $D_Q=D_P+d_{e}$\\
$\tau_P$ & \makecell[tl]{Earliest departure time assuming\\maximum speed} & $\max(e_i,d_e/b_k)+s_i$ & $\tau_Q=\max(e_j,\tau_P+d_e/b_k)+s_j$\\
$S_{P}$ & Total service time before $n_P$ & 0 & $S_{Q}=S_{P}+s_{n_{P}}$\\
${\cal T}_P$ & \makecell[tl]{Triples of distance, total service time\\and earliest service time\\for calculating $\sigma_P$} & $\{ (d_e, 0, e_i) \}$ & ${\cal T}_Q = {\cal T}_P \cup \{ (D_Q, S_Q, e_j) \}$\\
$\sigma_P$ & \makecell[tl]{Minimum vehicle speed such that\\all time windows are respected} & $d_e/l_i$ & $\sigma_Q = \max\limits_{(d,s,t) \in {\cal T}_P} \frac{D_Q-d}{l_{j}-t-(S_Q-s)}$\\
$\nu_P$ & Optimal speed & $\max(v^*_k,d_e/l_i)$ & $\nu_Q = \max(v^*_k,\sigma_Q) $\\\makecell[tl]{
$\alpha_P$\\$\beta_P$\\$\gamma_P$\\$\delta_P$} & \multirow{4}{*}{Coefficients for $\text{CO}_2e$ emissions} & \makecell[tl]{$\alpha_{ke}$\\$\beta_{ke}$\\$\gamma_{ke}$\\$\beta_P q_i$} & \makecell[tl]{$\alpha_Q = \alpha_P + \alpha_{ke}$\\$\beta_Q = \beta_P + \beta_{ke}$\\$\gamma_Q = \gamma_P + \gamma_{ke}$\\ $\delta_Q=\delta_P+\beta_Q q_j$}\\
\bottomrule
\end{tabular}

\footnotesize{$^a$ Assuming label $P$; $^b$ Assuming vehicle $k$ and a label representing the partial route along arc $e=(0,i)$; $^c$ Assuming label $Q$ obtained by extending $P$ along arc $e=(i,j)$.}
\end{minipage}
\end{table}

A key component of our pricing algorithm is the set of label extension procedures that do not require full backtracking to determine the cost of a route under the optimal speed. The idea is formalized with the following definition and proposition.

\begin{definition}[Cost of a path]
Let $Q$ be a path with arcs $(e_1, e_2, ..., e_L)$ and nodes $(v_0, v_1, v_2, ..., v_L)$ where $v_0=v_L=0$ is the depot. The \emph{cost of path} $Q$ when travelling at a speed of $\nu$ is defined as
\begin{align}
 c_k\Big(\frac{\alpha_Q }{\nu}  +  \delta_Q  + \beta_Q w_k  + \gamma_Q \nu^2\Big)  - \pi_Q, \label{path cost 1}
\end{align}
where $\alpha_Q=\sum_{l=1}^L \alpha_{ke_l}$, $\beta_Q =\sum_{l=1}^L \beta_{ke_l}$,  $\gamma_Q=\sum_{l=1}^L \gamma_{ke_l}$,  $\pi_Q=\sum_{l=1}^L \pi_{e_l}$, $\delta_Q =\sum_{i=1}^{L-1} \beta_{P_i} q_i$, and $P_i$ denote the path $(v_0, v_1, v_2, ..., v_i)$.
\end{definition}

\begin{proposition}
Let $P$ be a complete path that ends at the depot, and let $p$ be the route induced by $P$. Then, the cost of $p$ executed under the optimal vehicle speed, given by $c_p$, is equal to the cost of $P$ as determined by \eqref{path cost 1} where $\nu=\nu_P$.
\end{proposition}
\proof{}
Let $Q$ be a partial path of vehicle $k$ with arcs $(e_1, e_2, ..., e_L)$ and nodes $(v_0, v_1, v_2, ..., v_L)$ where $v_0=v_L=0$ is the depot with demand $q_0$ set to $0$. The travel cost of the route induced by $Q$, according to the objective function \eqref{sub: obj}, is given by
\begin{align}
c_k\sum_{l=1}^L\Big(\frac{\alpha_{ke_l}}{\nu}+ \beta_{ke_l}(w_{k}+x_{e_l})+\gamma_{ke_l}\nu^2\Big)-\sum\limits_{l=1}^L \pi_{v_l}, \label{pf: routing cost1}
\end{align}
where $x_{e_l}=\sum_{m=l}^{L-1} q_m$ is the sum of the demand of the customers in the remaining part of the route, that is, the payload of arc $e_l$. As shown below, the cost calculation by \eqref{pf: routing cost1} is equivalent to \eqref{path cost 1}.

The travel cost \eqref{pf: routing cost1} can be rewritten as
\begin{align*}
&c_k \Big(   \frac{ \sum_{l=1}^L\alpha_{ke_l}}{\nu}  + \sum\limits_{l=1}^L\sum_{m=l}^L \beta_{ke_l} q_m+ \sum\limits_{l=1}^L\beta_{ke_l} w_k + \sum\limits_{l=1}^L \gamma_{ke_l} \nu^2\Big) - \sum\limits_{l=1}^L \pi_{v_l} \\
=\ &c_k\Big(\frac{\alpha_Q }{\nu}  + \sum\limits_{l=1}^L\sum\limits_{m=l}^L \beta_{ke_l} q_m+ \beta_Q w_k  + \gamma_Q \nu^2\Big)  - \pi_Q.
\end{align*}

Note that $\sum_{l=1}^L\sum_{m=l}^L\beta_{ke_l}q_m=\sum_{l=1}^L\sum_{m=1}^l\beta_{ke_m}q_l=\sum_{l=1}^L \beta_{P_l}q_l=\delta_Q$, and hence equivalent to the routing costs defined in the objective function \eqref{sub: obj}.
\Halmos
\endproof

The pricing problem can be considered as a variant of the resource-constrained shortest-path (RCSP) problem \citep[see e.g.][]{feillet2004exact}. Typically, one finds the RCSP with a labeling procedure where dominance rules are employed to discard non-promising partial paths. In our case, however, dominance rules are likely to be ineffective because of the large number of resources involved. More specifically, in addition to the usual resources to handle vehicle capacity and time windows, in our case one must also observe dominance conditions on the allowed vehicle speed range and each cost component individually. Therefore, in our pricing algorithm, we decide to control the combinatorial growth of labels exclusively with completion bounds, which are detailed next.

\subsubsection{Completion bounds}
A completion bound is a lower bound on the reduced cost of all routes that can be generated from a label. Completion bounds accelerate the solution of the pricing problem since partial paths with nonnegative bounds are discarded during the labeling procedure.

Consider a partial path $P$ with arcs $(e_1, e_2, ..., e_L)$ and nodes $(v_0, v_1, v_2, ..., v_L)$. The precise cost along $P$ depends on the customers visited after $v_L$, since the cost along each arc depends on the arc payload. A lower bound on the cost along $P$, however, can be computed as follows:
\begin{equation}\label{eq:lbcostpp}
    \Phi_{P}=f_{k}+c_k\sum_{l=1}^{L}\Big(\frac{\alpha_{ke_l}}{\nu}+ \beta_{ke_l}(w_{k}+\bar{x}_{e_l})+\gamma_{ke_l}\nu^2\Big),
\end{equation}
where
\begin{equation}\label{eq:payload}
    \bar{x}_{e_l}=
    \begin{cases}
    q_P-\sum_{m=1}^{l-1}q_{v_{m}},&\text{if }\beta_{ke_l}\geq 0,\\
    Q_k-\sum_{m=1}^{l-1}q_{v_{m}},&\text{if }\beta_{ke_l}<0.
    \end{cases}
\end{equation}

Equation \eqref{eq:payload} considers a best-case (i.e., cost-minimizing) scenario concerning the payload along each arc. If the corresponding $\beta$ is nonnegative, the vehicle is assumed to travel along arc $e$ as light as possible. Otherwise, the vehicle is assumed to travel as loaded as possible.

Given the lower bound \eqref{eq:lbcostpp}, we propose two completion bounds for a label $P$. The first bound is based on a RCSP and explores the capacity resource to find a lower bound on the reduced cost of any extension of $P$. The second bound is based on a knapsack problem and explores not only the capacity resource but also the ``timing'' resources, that is, the fact that customers cannot be visited after their time windows and the vehicle must return to the depot no later than instant $l_{0}$.

We start with the RCSP-based completion bound, which adapts the bound proposed by \cite{florio2021routing} for solving the elementary RCSP. First, we associate to each arc $e=(i,j)\in\cal A$ a lower bound $\bar{\phi}_{e}$ on the reduced cost change when partial path $P$ is extended along $e=(i,j)$:
\begin{equation}\label{eq:rcspcosts}
    \bar{\phi}_{e}=\frac{\alpha_{ke}}{v^*_k}+ \beta_{ke}(w_{k}+x'_{e})+\gamma_{ke}(v^*_k)^2-\pi_{j},
\end{equation}
where $x'_{e}=0$ if $\beta_{e}\geq 0$ and $x'_{e}=Q_{k}$ otherwise.

We denote by $S^{*}_{i}(Q)$ the value of the RCSP from node $i\in\cal N\setminus\{0\}$ to node 0 in a graph with arc costs given by \eqref{eq:rcspcosts}, in which the initial resource limit is $Q$ and an amount $q_{j}$ of resource is consumed each time node $j\neq0$ is visited. Then, the RCSP-based completion bound is given by:
\begin{equation}\label{eq:rcspbound}
    \Phi_{P}-\sum_{i\in N_{P}}\pi_{i}+S^{*}_{n_{P}}(Q_{k}-q_{P}).
\end{equation}

Equation \eqref{eq:rcspbound} yields a valid bound because $\Phi_{P}-\sum_{i\in N_{P}}\pi_{i}$ is a lower bound on the reduced cost of partial path $P$, and $S^{*}_{n_{P}}(Q_{k}-q_{P})$ is a lower bound on the reduced cost of any feasible extension to $P$. To enable evaluations of \eqref{eq:rcspbound} in constant time for any label $P$, at the beginning of an iteration of the pricing problem we pre-compute $S^{*}_{i}(Q)$ for all $i\in\cal N\setminus\{0\}$ and $Q\in\{0,\ldots,Q_{k}\}$. The (non-elementary) RCSPs can be solved efficiently by dynamic programming.

While the RCSP bound is computationally efficient, it does not explore time window constraints nor the fact that we price elementary routes. In the knapsack bound, we set up a $\{0,1\}$-knapsack problem with capacity of $l_{0}-\tau_{P}$, which corresponds to the maximum remaining routing time after customer $n_{P}$ is served. Then, we define a set of knapsack items
\begin{equation}\nonumber
    \mathcal{I}=\{i\in\mathcal{N}\setminus(\{0\}\cup N_{P}):q_{i}+q_{P}\leq Q_{k}\wedge\tau_{P}+d_{n_{P},i}/b_{k}\leq l_{i}\}.
\end{equation}

Each element of $\mathcal{I}$ is a customer that can be visited after $n_{P}$, considering time window and vehicle capacity constraints. With each item $i\in\mathcal{I}$ a value $v(i)$ and a weight $w(i)$ are associated:
\begin{align}
    v(i)&=\max_{e=(j,i)\in\cal A}-\bar{\phi}_{e},\nonumber\\
    w(i)&=\min_{(j,i)\in\mathcal{A}}d_{ji}/b_{k}.\nonumber
\end{align}

The value and weight of an item $i$ correspond to the maximum reduced cost decrease and minimum amount of the time resource consumed, respectively, when customer $i$ is visited in an extension of partial path $P$. We let $K^{*}_P$ be the optimal solution value of the knapsack problem defined above. Then, the knapsack completion bound is given by:
\begin{equation}\label{eq:kpbound}
\Phi_{P}-\sum_{i\in N_{P}}\pi_{i}-K^{*}_P.
\end{equation}

Each time a label $P$ is generated, we evaluate \eqref{eq:kpbound} and discard $P$ if the bound is nonnegative. This evaluation requires solving the knapsack problem to optimality, which can also be achieved efficiently by dynamic programming.

\subsection{Branch-and-bound}
The implemented branch-and-bound framework finds an optimal integer solution to (SP) by branching on variables $y_e$, $e\in\mathcal{A}$, that take fractional values. We apply a semi-strong branching rule where each potential branching variable is evaluated under the current pool of columns, and the variable on which branching leads to the highest lower bound is chosen. More precisely, at a given branch-and-bound node, we let $\Omega_{\textsf{R}}$ be the set of all columns generated and $\cal Y$ the set of arc variables that assume fractional values in the solution to the RMP. Then, we evaluate
\begin{equation}\label{eq:evalbranching}
\min\{\text{RMP}(\Omega_{\textsf{R}},\{y_{e}\}, \{\}),\text{RMP}(\Omega_{\textsf{R}},\{\},\{y_{e}\})\}
\end{equation}
for each $y_{e}\in\cal Y$, where RMP$(\Omega,\cal Y_0,\cal Y_1)$ corresponds to the optimal solution value of the RMP restricted to columns $\Omega$ and enforcing constraints $y_e=0$ for all $y_e\in\cal Y_0$ and $y_e=1$ for all $y_e\in\cal Y_1$ in addition to the branching constraints of the parent branch-and-bound node. Finally, we branch on the variable $y_e\in\cal Y$ such that \eqref{eq:evalbranching} is maximum. Note that evaluating \eqref{eq:evalbranching} for each potential branching variable can be implemented efficiently by loading a single linear program and (re)solving it for each variable after adjusting the cost vector accordingly, by penalizing the cost of routes that do not comply with the candidate branching decision.

\section{Computational Experiments}\label{sec:experimental-results}

In this section, we will evaluate the performance of the tabu search heuristic described in Section  \ref{sec: metaheuristics} (denoted as \textbf{TS}) and the branch-and-price approach described in Section \ref{sec: BCP} (denoted as \textbf{BP}).
Afterward, we make use of the heuristic in an empirical study for evaluating the potential benefits of using elevation data in optimizing the vehicle routes.

We organize the remainder subsections as follows.
Section \ref{sec: instances} describe the test instances constructed using data from literature, and Section \ref{subsec:test-instances-using-real-world-data} describe the test instances constructed using real-world data.
Section \ref{subsec:-parameter-tuning} is about the parameter tuning experiment.
In Section  \ref{sec: results PRP-SO}, we evaluate the efficiency of BP and TS.
In Section \ref{subsec:the-value-of-using-elevation-information}, we evaluate the potential benefits of using elevation data for planning the vehicle routes.

\subsection{Test instances using data from literature} \label{sec: instances}
Since there are no existing benchmark datasets of PRP-SO, we constructed test instances based on the instances of \cite{solomon1987algorithms} which are originally created for VRPTW.
The VRPTW instances are adapted into PRP-SO instances by associating randomly generated elevation information on the nodes.
The units for time and demand are also scaled to match the realistic instances.

The VRPTW instances have four sets of instances involving 25, 50, 75, and 100 customers, respectively.
Instances are divided into three classes according to the customer location distribution: clustered distribution (C class), scattered distribution (R class), and partially scattered and partially clustered distribution (RC class).
Each class is further subdivided into the narrower time window class and the wider time window class.
In our experiments, we will test on the instances with narrower time windows.
In total, there are 116 instances of PRP-SO constructed from using the instances of \cite{solomon1987algorithms}.
The remaining part of this subsection describes how the VRPTW instances are adapted into the PRP-SO instances.

The elevation of the nodes is randomly generated with a uniform distribution between 0 and 1,000 meters.
The distance in kilometers between node $i$ and node $j$ is given by $$D_{ij} = \sqrt{(X_i - X_j)^2 + (Y_i-Y_j)^2 + \big(\frac{Z_i - Z_j}{1000}\big)^2}$$ which rounds to the nearest meter, where $(X_i, Y_i)$  and $(X_j, Y_j)$ are the coordinates, and $Z_i$ and $Z_j$ are the elevations of node $i$ and $j$ respectively.
With the rounded values of distances, the road angle between two nodes is given by
$\tan^{-1}\Big( \frac{Z_i - Z_j}{ D_{ij} }\Big). $

The units for time and demand are also scaled to match the realistic instances.
For our experiments, service times and the time windows have a unit of 0.02 hours.
For example, a due date of 1236 from the Solomon's instances represents a due date at the 24.72 hours (given by 1236*0.02) after the planning horizon starts.
This implies that, if vehicles always travel at 50 km per hour, the time window constraints would remain the same as the ones in the original VRPTW instances.

\begin{table}[!ht]
 \caption{Vehicle Information} \label{data vehicles}\small\centering
 \begin{tabular}{cccccc}
  \toprule
Original capacity &  Vehicle type & Capacity & Demand unit   \\
  \midrule
200 &  LDV& 1,200 kg&6 kg per unit          \\
700 &  MDV &12,600 kg &18 kg per unit        \\
1,000 &  HDV &31,000 kg &31 kg per unit        \\
  \bottomrule
 \end{tabular}
\end{table}

Three types of vehicles appeared in Solomon's instances: 200, 700, and 1,000 units.
We scale the demand and the vehicle capacity accordingly so that it is equivalent to the original constraints for vehicle capacity and at the same time matches typical truck classifications: LDV, MDV, and HDV  for the 200, 700, and 1,000 capacity units respectively.
Table \ref{data vehicles} summarizes the vehicle capacity and demand unit used in our experiments.
For example, the vehicles in the Solomon’s instances with a capacity of 200 units correspond to the LDV vehicles of the PRP-SO instances, and therefore a demand of 20 units in those instances represents a demand of 120 kg (given by $20\times 6$) in the PRP-SO instances.

All the vehicles have a fixed cost of 100 EUR, a fuel cost of 1.42 EUR per liter, and a maximum speed of 80 km per hour.
Other parameters used for the $\text{CO}_2e$ consumption calculations are summarized in Table  \ref{CMEM parameters}.

\subsection{Test instances using real-world data}\label{subsec:test-instances-using-real-world-data}
\begin{figure}[h!] 
\FIGURE{
\includegraphics[width=10cm, height=6cm]{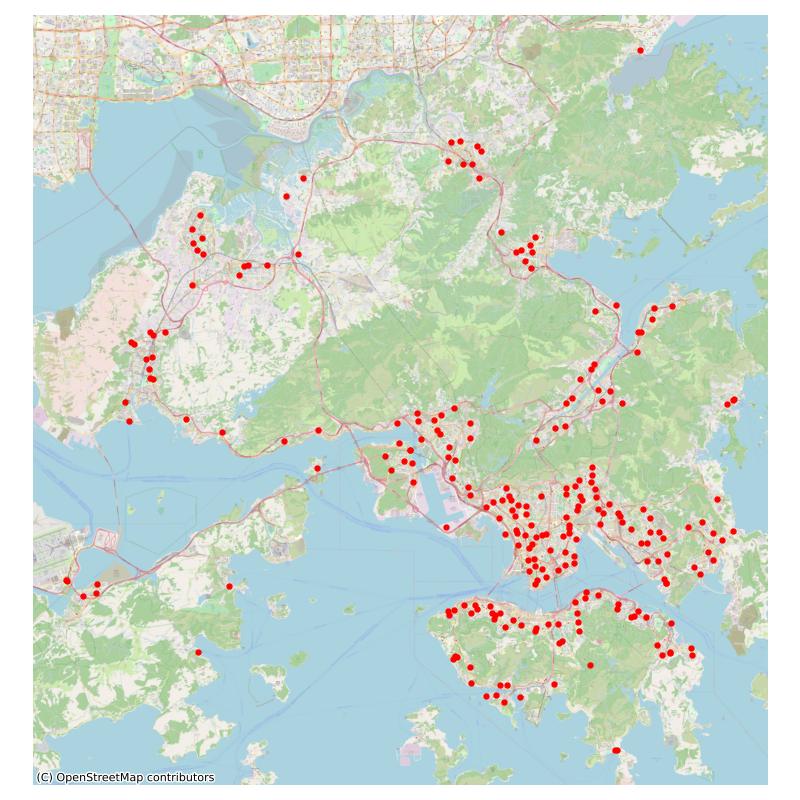}
}
{Customer locations \label{fig: hkmap}}
{}
\end{figure}

Another dataset is constructed based on the distribution network of a large international health and beauty retailer.
There are 248 customers considered in our experiments which are representing the retailer's stores in Hong Kong.
Geometric information is obtained by using Google Maps APIs, including the coordinates, elevations, suggested paths between the stores.
Figure \ref{fig: hkmap} illustrates the geographical locations of these stores.
The landscape varies from fairly hilly to mountainous with steep slopes.
The stores' locations are clustered, densely populated in the central areas.
We will use this dataset for evaluating the potential benefits of using elevation information in optimizing the vehicle routes.

We preprocess the geographic data from the Google Maps API and the Elevation API into the $\alpha$, $\beta$, and $\gamma$ values associated on the arcs (defined in Section \ref{sec:the-hill-climbing-vehicle-routing-problem}), so that the proposed solution approaches can compute the fuel consumption and  $\text{CO}_2e$ emission efficiently.
To begin with, for every distinct pair of the stores, we obtain a suggested path by using the Google Maps API and find out the elevation for all the coordinates along the suggested path by using Google Elevation API\@.
Afterward, to construct arc segments, coordinates along a path are divided into segments, with each segment the distance is no longer than 1,000 meters.
An arc segment can be viewed as a slope along the suggested path.
In our experiments, there are in total 155,333 such slopes.
Lastly, the angles and distance of all these slopes are determined using the coordinates and elevation data, and thus we have the $\alpha$, $\beta$, and $\gamma$ values of all the arcs connecting the stores.

In our experiments, nine instances are constructed.
Each instance consists of 100 customers which are randomly selected amongst the 248 stores.
The ready time, due time, demand, vehicle number, and capacity are data from the Solomon's C class instances.
The depot is located at the Kwai Tsing Container Terminal which is the busiest port in Hong Kong.

\subsection{Parameter tuning} \label{subsec:-parameter-tuning}
The best parameters amongst the values specified in Table \ref{tab: parameter values} are chosen for each instance class.

\begin{table}[!ht]\caption{Parameter Values}  \label{tab: parameter values}
\centering
\begin{tabular}{cllc}
 \toprule
 Parameter  &Possible values& Description \\
 \midrule
 $I_1$ & 1, 2 & Number of iterations in the first phase\\
 $I_2$ & 100,000 & Number of iterations in the second phase  \\
 $\lambda$ & $10^{-6}$, $5\times 10^{-6}$, $10^{-7}$, $5\times 10^{-7}$ & Diversification intensity  \\
 $h$ &3, 4, 5, 6& Parameter for setting the tabu tenure \\
 $\delta$  & 10, 20, 30, 40 &Penalty update frequency \\
 \bottomrule
\end{tabular}
\end{table}

All experiments have been conducted on a server computer running Ubuntu with an Intel Xeon CPU E5-2698 v3 @ 2.30GHz, 16 cores, and 15 GB of main memory.
Algorithms have been implemented in C++ and compiled using GNU g++ version 10.2.0 with -O2 flag.
The algorithms run on a single core per instance.

\subsection{Efficiency of the approaches}  \label{sec: results PRP-SO}
Table \ref{HC-VRPTW results} summarizes the computational results of BP and TS on the instances described in Section \ref{sec: instances} with 25, 50, 75 and 100 customers respectively. 
Results on individual instances are shown in Appendix \ref{app: results}.
Column \textbf{NI} is the number of instances in the dataset class excluding the ones that no feasible solution can be found by using BP within 3 hours.
Column \textbf{NO} is the number of instances that can be solved to optimality by using BP within the time limit.
Column \textbf{NV} is average number of vehicle routes, column \textbf{TD} is average total travel cost, column \textbf{AT} is the average cpu time (in seconds), and  column \textbf{AG} is the average optimality gap (in percentage).
When reporting the average values, we excluded the instances that no feasible solution can be found by using BP within 3 hours.
\begin{table}[!ht] \small\centering
\caption{Computational Results of the PRP-SO Instances \label{HC-VRPTW results}}
\begin{tabular}{llc|rrrrr|rrrrrrrrrrrrrrrrr}
 \toprule
 & & &\multicolumn{4}{c}{BP}& & \multicolumn{3}{c}{TS}  \\
 \multicolumn{2}{c}{Class}   &NO/NI  &NV & TD & AT (s) & AG (\%) && NV & TD &  AT (s) \\
 \midrule
 N25 & C  & 9/9  & 3.000 & 24.059 & 215.90 & 0\% & & 3.000 & 24.597 & 0.41 \\
 & RC & 8/8  & 3.250 & 45.314 & 30.28  & 0\% & & 3.250 & 45.331 & 0.98 \\
 & R  & 12/12 & 4.667 & 60.240 & 30.02  & 0\%  && 4.750 & 59.306 & 2.49 \\
 N50 & C  & 7/7  & 5.000 & 45.641  & 3608.80 & 0.000\%  && 5.000 & 46.612  & 2.65  \\
 & RC & 4/8  & 6.500 & 94.850  & 6933.27 & 0.087\% & & 6.500 & 97.514  & 10.26 \\
 & R  & 9/10 & 7.500 & 107.889 & 2707.29 & 0.018\% & & 8.000 & 109.131 & 8.17 \\
 N75 & C  & 2/4  & 8.000  & 83.020  & 8314.54 & 0.079\% & & 8.000  & 82.994  & 3.85  \\
 & RC & 1/6  & 9.667  & 154.001 & 9262.17 & 4.088\% & & 10.333 & 153.598 & 10.03 \\
 & R & 5/6  & 11.833 & 152.447 & 3956.78 & 0.018\% & & 12.500 & 160.035 & 8.23 \\
 N100 &C  & 1/1 & 10.000 & 107.742 & 3255.73 & 0.000\% & & 10.000 & 107.742 & 2.21  \\
 &RC & 1/4 & 12.750 & 210.165 & 8957.57 & 2.440\% & & 13.000 & 195.959 & 12.93 \\
 &R  & 2/3 & 16.667 & 194.144 & 4489.88 & 1.983\% & & 13.500 & 170.702 & 42.31 \\
 \midrule
Total: & & 61/78 & 98.8 & 1279.5 & 51762.2 & & &97.8 & 1253.5 & 104.5 \\
\bottomrule  
\end{tabular}
\end{table}

As shown in Table \ref{HC-VRPTW results}, BP can solve instances up to 100 customers, with 61 instances (53\%) solved to optimality.
BP can find all the optimal solutions of the dataset with 25 customers within a reasonable time, and most optimal solutions of the dataset with 50 customers within 3 hours.
As compared to TS, BP can determine better solutions on smaller instances, but at the expense of significantly more CPU time.
TS can find near-optimal solutions for all instances within one minute, and outperforms BP on solution quality for the dataset with 100 customers.

%

\subsection{The value of using elevation information}\label{subsec:the-value-of-using-elevation-information}
\begin{figure}[!ht]
\FIGURE{
 \includegraphics[height=5cm]{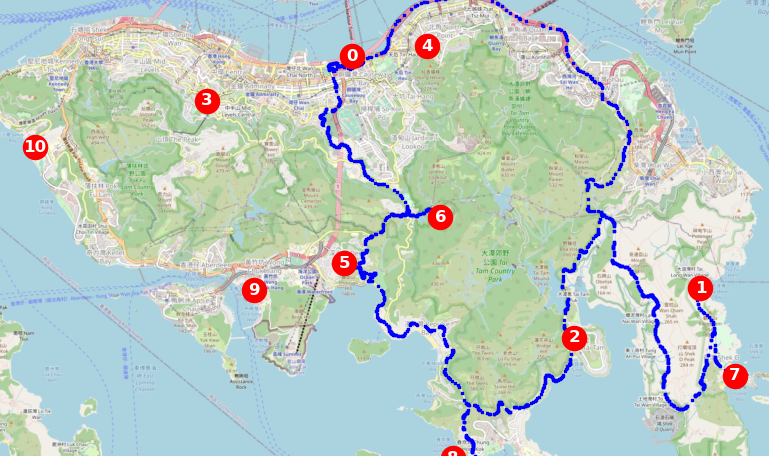} \ \  \includegraphics[height=5cm]{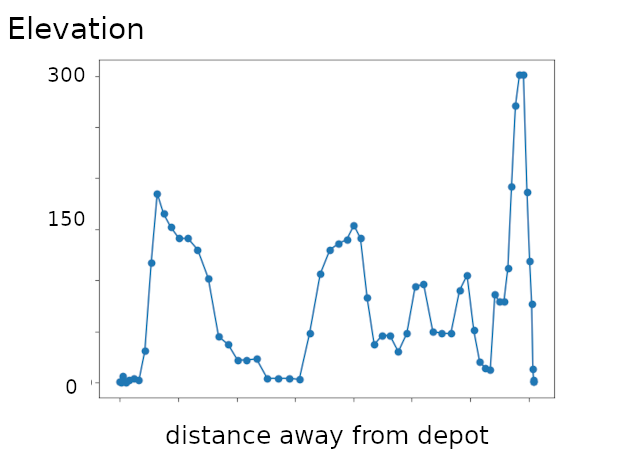}
}{An example vehicle route \label{fig: Example solutions}}{}
\end{figure}

For evaluating the potential benefits of using elevation information in optimizing the vehicle routes, the real-world instances described in Section \ref{subsec:test-instances-using-real-world-data} are solved by using TS.
Table \ref{tab: dataset B results} summaries the total distance (in km), average speed (in km/h), total fixed cost (in EUR), total fuel cost (in EUR), and the total elevation (in meters).
Figure \ref{fig: Example solutions} illustrates an example vehicle route.
The same dataset is solved again with which elevation information is ignored when planning the vehicle routes.
The results are reported in Table \ref{tab: dataset B results ignore slope}.
This is achieved in our experiment by setting the angles of all slopes to zero when optimizing the vehicle routes by using TS\@, and evaluating the solutions with the correct elevations and slopes after the vehicle routes have been decided.

\begin{table}[!ht] \caption{Results on Real-world Instances} \label{tab: dataset B results}
\centering
\begin{tabular}{cccccc}
 \toprule
 {Instance}& {Distance}& {Speed}& {Fixed cost}& {Fuel cost}& {Elevation} \\
 \midrule
 HK01              & 1638.6                       & 62.3                          & 1500                      & 109.5                    & 1816                               \\
 HK02              & 1460.1                       & 62.3                          & 1400                      & 163.6                    & 1436                               \\
 HK03              & 990.2                        & 61.1                          & 1200                      & 196.7                    & 1070                               \\
 HK04              & 1153.7                       & 63.1                          & 1000                      & 50.0                     & 1440                               \\
 HK05              & 1209.0                       & 63.5                          & 1400                      & 197.5                    & 1185                               \\
 HK06              & 1441.9                       & 64.9                          & 1300                      & 89.4                     & 1748                               \\
 HK07              & 1159.5                       & 65.5                          & 1300                      & 185.4                    & 1309                               \\
 HK08              & 1384.8                       & 64.6                          & 1200                      & 110.3                    & 1359                               \\
 HK09              & 1088.5                       & 63.7                          & 1100                      & 81.5                     & 1495                               \\
 \midrule
Average:           & 1280.7                       & 63.5                          & 1266.7                    & 131.5                    & 1428.7\\
 \bottomrule
\end{tabular}
\end{table}

\begin{table}[!ht]\caption{Results on Real-world Instances when Slopes are Ignored} \label{tab: dataset B results ignore slope}
\centering
 \begin{tabular}{cccccc}
  \toprule
  {Instance}& {Distance}& {Speed}& {Fixed cost}& {Fuel cost}& {Elevation} \\
  \midrule
  HK01    & 1158.7 & 60.6 & 1500   & 405.1 & 1582   \\
  HK02    & 1015.1 & 60.7 & 1400   & 331.1 & 1152   \\
  HK03    & 953.0   & 60.9 & 1200   & 184.5 & 1043   \\
  HK04    & 765.2  & 60.5 & 1100   & 114.5 & 1443   \\
  HK05    & 1040.7 & 62.9 & 1400   & 409.8 & 1140   \\
  HK06    & 1071.8 & 62.8 & 1300   & 518.3 & 1536   \\
  HK07    & 1004.4 & 62.3 & 1300   & 278.1 & 1311   \\
  HK08    & 945.3  & 64.0 & 1300   & 197.7 & 1307   \\
  HK09    & 831.7    & 59.5 & 1100   & 138.0 & 1331   \\
  \midrule
  Average:           & 976.2    & 61.6 & 1288.9 & 286.3  & 1316.1\\
  \bottomrule
 \end{tabular}
\end{table}


We can observe the impacts on solutions when elevation information is used for planning the vehicle routes by comparing the results in Table \ref{tab: dataset B results} and Table \ref{tab: dataset B results ignore slope}.
As shown in the experimental results, when elevation information is considered in planning, fuel consumption decreased by 54\% on average while the average vehicle speed and elevation increased slightly, and the total travel distance increased by 31\%.
Larger travel distances and lower fuel consumption seem to be contradictory for typical vehicle routing problems.
This is because differs from typical vehicle routing problems, the fuel consumption now depends not only on the distance, but also on the slopes, payload, and vehicle speeds.
Optimal vehicle routes should therefore save fuel costs by avoiding going uphill at a high speed with a large payload.
As a result, vehicles tend to visit more customers first before going uphill.
Although this will increase travel distance, fuel consumption can be saved from going uphill with less payload.

For typical vehicle routing problems or when vehicle routes are planned manually, travel distance is usually minimized in the objective function.
Our experimental result reveals that this can lead to a suboptimal solution in reality.
With the elevation information, the optimal solution can balance the tradeoff between the energy consumption due to longer distances and the higher payload when going uphill.
To avoid high fuel consumption when vehicles going uphill, heavier items tend to be delivered first before going uphill.
Customer time windows have to be taken into account so that vehicles do not need to speed up to meet the due times which can result in high fuel consumption.
If elevation information is ignored when planning the vehicle routes, fuel consumption due to payload is underestimated when vehicles are going uphill, which leads to poor solutions.
With the significant savings we observed from the experimental results, logistic service providers should consider using elevation data for planning their vehicle routes in practice.

\subsection{Impact of payloads and slopes}
\begin{figure}[!ht]
\FIGURE
{\includegraphics[height=5.3cm]{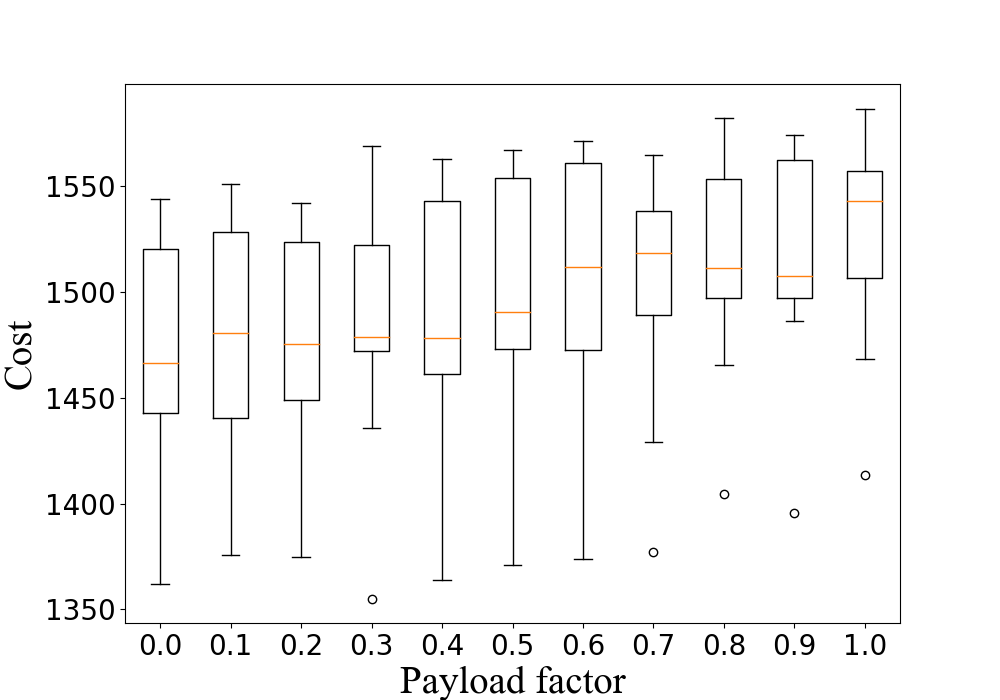} \ \  \includegraphics[height=5.3cm]{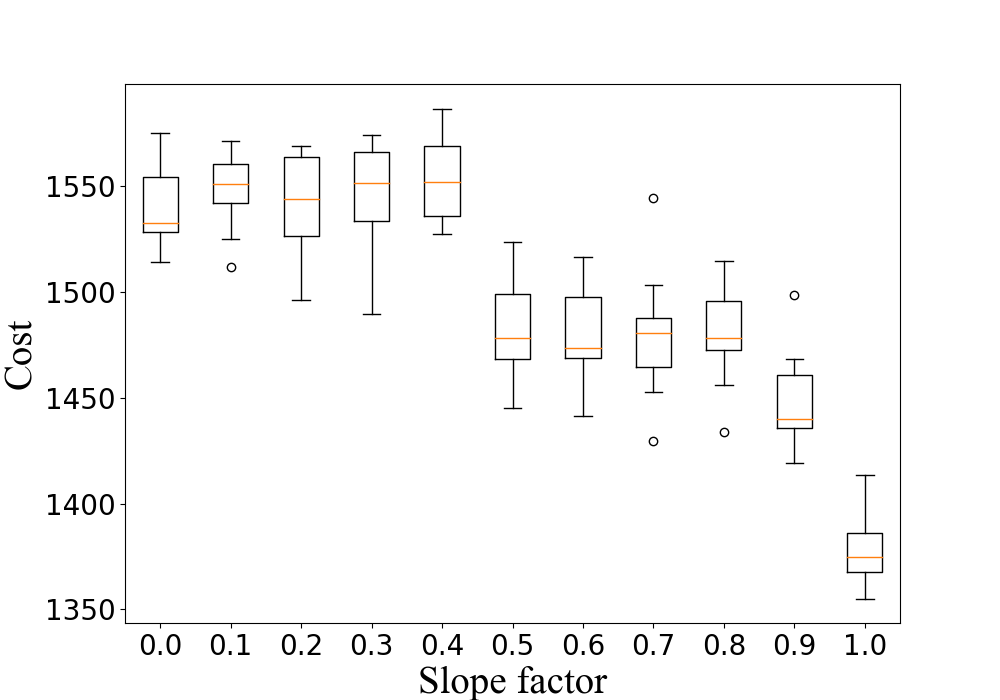}
}
{Impact of payloads and slopes on routing costs \label{fig:  impact of payloads and slopes on the routing costs}}
{}
\end{figure}

For evaluating the impact of payloads and slopes on the optimized solutions, the real-world instances described in Section \ref{subsec:test-instances-using-real-world-data} are modified by scaling the payloads by a factor $r_1 \in \{0, 0.1, ..., 1\}$ when calculating the costs, and scaling the slopes by a factor $r_2 \in \{0, 0.1, ..., 1\}$.
In our experiments,  we replace the payload $x_{ke}$ in \eqref{co2 cost} by a factor $r_1 x_{ke}$ when calculating the costs, and replace the slope $\phi_{es}$ in \eqref{co2 beta} by  $r_2\phi_{es}$ when preprocessing the data.
A higher value of the payload factor ($r_1$) represents scenarios when relatively heavier items are shipped, and a higher value of the slope factor ($r_2$) represents more hilly areas where steep slopes commonly appear.
There are 1089 instances tested in total, and each of them is solved by using TS with a timelimit of 300 seconds.
Table \ref{tab: impact of payloads and slopes} summarizes the average costs (and refer to appendix \ref{tab: impact of payloads and slopes full} for the complete results).
Figure \ref{fig:  impact of payloads and slopes on the routing costs} shows  the average cost with varying payloads ($r_1$) and slopes ($r_2$) respectively. 
	
\begin{table}[!ht]\caption{Impact of Payloads and Slopes on Routing Costs} \label{tab: impact of payloads and slopes}
	\centering 
\begin{tabular}{cc|ccccccccccc|c}
                              &         & \multicolumn{7}{c}{Slope factor ($r_2$)}                                \\ \hline
                              &         & 0      & 0.2    & 0.4    & 0.6    & 0.8    & 1      & Average \\ \hline
\multirow{7}{*}{\makecell[l]{Payload\\factor\\ ($r_1$)}} & 0       & 1516.0 & 1543.9 & 1527.4 & 1458.7 & 1456.2 & 1361.8 & 1477.3  \\
                              & 0.2     & 1514.3 & 1519.5 & 1527.3 & 1472.0 & 1475.4 & 1374.8 & 1480.6  \\
                              & 0.4     & 1530.7 & 1563.0 & 1540.4 & 1466.4 & 1478.4 & 1363.9 & 1490.5  \\
                              & 0.6     & 1545.9 & 1565.5 & 1556.6 & 1475.9 & 1512.0 & 1373.8 & 1505.0  \\
                              & 0.8     & 1553.0 & 1532.5 & 1582.0 & 1511.1 & 1494.6 & 1404.7 & 1513.0  \\
                              & 1       & 1575.1 & 1542.8 & 1586.7 & 1516.6 & 1496.1 & 1413.3 & 1521.8  \\ \hline
                              & Average & 1539.2 & 1544.6 & 1553.4 & 1483.5 & 1485.5 & 1382.0 &        
\end{tabular}
\end{table}

From the experimental results, we can see that the shipping cost increases with the payload ($r_1$) and decreases with the slope ($r_2$).
It is more costly to ship with a heavier payload (higher value of $r_1$) with a percentage increase in total costs up to 3.49\% on average regardless of the slopes.
It is less costly to ship in more hilly areas (higher value of $r_2$) with a saving up to 11.71\% of the total costs on average regardless of the payloads. 
The impact due to slopes is significantly higher than the impact due to payloads, which reveals the importance of taking into account slopes when optimizing the vehicle routes.

%
%

%

\section{Conclusions} \label{sec: conclusions}


This paper formulates and proposes efficient solution methods for a joint Pollution Routing and Speed Optimization Problem (PRP-SO), where the total travel cost is a function of fuel consumption and $\text{CO}_2e$ emissions and depends, simultaneously, on road grades, arc payloads, and vehicle speed. 
The introduction of the vehicle speed as continuous decision variables results in more complicated optimization subproblems in the presence of time window constraints. For the fixed-sequence speed optimization problem where the vehicle route is known, the proposed approach is conceptually simple and computes the optimal vehicle speed (with and without time windows) in quadratic time. In the speed optimization with variable routes, we introduced a novel labeling algorithm, without full backtracking searching, that efficiently determines the cost of a vehicle route that travels with optimal speed.
Based on the proposed speed optimization algorithms, we present two general solution approaches for solving the PRP-SO. An approximate solution strategy aims to solve large instances in a short computational time. For this purpose, we integrated the fixed-sequence speed optimization algorithm in a tabu search metaheuristic. The second approach consists of an exact branch-and-price algorithm, in which the variable-route speed optimization is managed within the pricing problem.

We carried out extensive computational experiments on modified Solomon benchmarks and newly constructed real-life instances. Numerical results show that the exact solution methodology performs very well in terms of solution quality: 61 out of 116 instances are solved to optimality. Contrary to the computational outcomes presented by \cite{dabia2017exact}, in which several instances with only 25 customers cannot be solved, we are able to solve all small-scale problem instances (25 customers) within a reasonable time. Our BP algorithm solved most of the problem instances with 50 customers and reached optimal solutions for some larger benchmark instances with up to 100 customers. We show that our metaheuristic works very effectively for all instances solved. The heuristic consumes less than one minute to find near-optimal solutions in all instances and improved best-known solutions where the exact algorithm did not reach optimality.

Our computational results on real-world instances provide sufficient evidence to suggest some essential managerial insights. First, significant savings (53\%) in fuel consumption and $\text{CO}_2e$ emissions are observed, especially when shipping heavy items in hilly areas. Second, vehicle routes included a larger number of customer visits (located at flatter terrain) before going to uphill destinations, also significantly reducing fuel consumption. Third, if elevation information is ignored when planning vehicle routes, fuel consumption estimation is inaccurate.

Finally, we do believe that pending issues are requiring future research. From a modeling perspective, the route security (uphill and downhill paths) in hilly topographic cities is a challenging issue to be included. Future studies can now be focused on generalizing the proposed methodological approaches to a new setting where two metrics of performance, fuel consumption, and vehicle route security should be optimized.    

    \bibliographystyle{informs2014trsc}
    \clearpage\newpage
    \bibliography{References}

    \begin{APPENDICES}{}
\clearpage\newpage
\section{Results} \label{app: results}
\begin{table}[!h]
\caption{Results on dataset A with 25 customers} \centering
\begin{tabular}{rrrrrrrrr}
\toprule
&\multicolumn{3}{c}{TS}  & \multicolumn{4}{c}{BP}        \\
Instance & NV & TD       & Time(s) & NV & TD       & Time(s) & Gap    \\
\midrule
c101  & 3 & 24.426 & 0.04  & 3 & 24.426 & 34.48  & 0.00\% \\
c102  & 3 & 23.836 & 0.33  & 3 & 23.278 & 167.03 & 0.00\% \\
c103  & 3 & 25.120 & 0.15  & 3 & 24.749 & 564.82 & 0.00\% \\
c104  & 3 & 25.835 & 0.15  & 3 & 23.994 & 750.72 & 0.00\% \\
c105  & 3 & 23.893 & 0.04  & 3 & 23.893 & 37.14  & 0.00\% \\
c106  & 3 & 22.743 & 0.04  & 3 & 22.743 & 30.45  & 0.00\% \\
c107  & 3 & 24.087 & 0.08  & 3 & 24.087 & 49.42  & 0.00\% \\
c108  & 3 & 25.666 & 1.10  & 3 & 25.199 & 86.84  & 0.00\% \\
c109  & 3 & 25.767 & 1.71  & 3 & 24.163 & 222.18 & 0.00\% \\
r101  & 8 & 75.787 & 0.17  & 8 & 75.787 & 3.43   & 0.00\% \\
r102  & 7 & 67.899 & 0.25  & 7 & 67.899 & 8.99   & 0.00\% \\
r103  & 4 & 60.948 & 25.69 & 4 & 60.948 & 20.41  & 0.00\% \\
r104  & 4 & 51.749 & 0.13  & 4 & 51.749 & 24.55  & 0.00\% \\
r105  & 5 & 73.139 & 2.24  & 5 & 71.878 & 12.81  & 0.00\% \\
r106  & 5 & 59.079 & 0.61  & 4 & 71.543 & 17.33  & 0.00\% \\
r107  & 4 & 52.313 & 0.08  & 4 & 52.313 & 39.68  & 0.00\% \\
r108  & 4 & 50.755 & 0.14  & 4 & 50.755 & 112.96 & 0.00\% \\
r109  & 4 & 58.009 & 0.13  & 4 & 58.009 & 17.49  & 0.00\% \\
r110  & 4 & 56.352 & 0.06  & 4 & 56.352 & 31.29  & 0.00\% \\
r111  & 4 & 55.261 & 0.14  & 4 & 55.261 & 29.14  & 0.00\% \\
r112  & 4 & 50.382 & 0.29  & 4 & 50.382 & 42.16  & 0.00\% \\
rc101 & 4 & 59.050 & 0.36  & 4 & 59.050 & 21.24  & 0.00\% \\
rc102 & 3 & 46.660 & 0.04  & 3 & 46.660 & 18.98  & 0.00\% \\
rc103 & 3 & 43.835 & 0.81  & 3 & 43.752 & 38.84  & 0.00\% \\
rc104 & 3 & 40.023 & 0.51  & 3 & 40.023 & 53.76  & 0.00\% \\
rc105 & 4 & 51.685 & 0.55  & 4 & 51.685 & 12.62  & 0.00\% \\
rc106 & 3 & 43.845 & 1.09  & 3 & 43.845 & 17.57  & 0.00\% \\
rc107 & 3 & 39.803 & 3.75  & 3 & 39.803 & 28.56  & 0.00\% \\
rc108 & 3 & 37.749 & 0.70  & 3 & 37.693 & 50.68  & 0.00\%\\
\bottomrule
\end{tabular}
\end{table}

\newpage
\begin{table}[!h]
\caption{Results on dataset A with 50 customers} \centering
\begin{tabular}{rrrrrrrrr}
\toprule
&\multicolumn{3}{c}{TS}  & \multicolumn{4}{c}{BP}        \\
Instance & NV & TD       & Time(s) & NV & TD       & Time(s) & Gap    \\
\midrule
c101  & 5  & 46.632  & 0.58  & 5  & 46.632  & 355.82   & 0.00\% \\
c102  & 5  & 44.457  & 1.41  & 5  & 44.324  & 9896.23  & 0.00\% \\
c105  & 5  & 45.762  & 3.34  & 5  & 45.762  & 1373.90  & 0.00\% \\
c106  & 5  & 43.440  & 3.14  & 5  & 43.440  & 1014.27  & 0.00\% \\
c107  & 5  & 45.942  & 1.45  & 5  & 45.942  & 1515.75  & 0.00\% \\
c108  & 5  & 49.758  & 5.09  & 5  & 48.022  & 4114.59  & 0.00\% \\
c109  & 5  & 50.292  & 3.51  & 5  & 45.368  & 6991.02  & 0.00\% \\
r101  & 11 & 153.622 & 40.54 & 11 & 136.449 & 66.09    & 0.00\% \\
r102  & 10 & 118.659 & 1.68  & 10 & 116.500 & 113.91   & 0.00\% \\
r103  & 8  & 109.075 & 9.20  & 8  & 100.962 & 352.12   & 0.00\% \\
r105  & 9  & 121.023 & 1.92  & 8  & 128.679 & 356.35   & 0.00\% \\
r106  & 8  & 105.068 & 5.83  & 7  & 110.509 & 646.84   & 0.00\% \\
r107  & 7  & 92.989  & 10.15 & 6  & 94.761  & 3521.55  & 0.00\% \\
r109  & 7  & 113.663 & 6.36  & 7  & 102.349 & 329.25   & 0.00\% \\
r110  & 7  & 94.271  & 3.59  & 6  & 102.664 & 2792.51  & 0.00\% \\
r111  & 7  & 95.106  & 2.04  & 6  & 100.870 & 8357.22  & 0.00\% \\
r112  & 6  & 87.839  & 0.36  & 6  & 85.143  & 10537.02 & 0.18\% \\
rc101 & 8  & 122.350 & 8.16  & 8  & 120.475 & 2761.42  & 0.00\% \\
rc102 & 7  & 112.088 & 21.62 & 7  & 108.434 & 10785.42 & 0.27\% \\
rc103 & 6  & 95.953  & 9.52  & 6  & 92.720  & 10772.43 & 0.22\% \\
rc104 & 5  & 73.271  & 19.88 & 5  & 71.435  & 2075.18  & 0.00\% \\
rc105 & 8  & 110.736 & 8.75  & 8  & 108.016 & 5045.61  & 0.00\% \\
rc106 & 6  & 96.262  & 3.35  & 6  & 94.017  & 10794.97 & 0.17\% \\
rc107 & 6  & 89.167  & 3.20  & 6  & 85.984  & 10736.23 & 0.04\% \\
rc108 & 6  & 80.282  & 7.58  & 6  & 77.722  & 2494.94  & 0.00\%\\
\bottomrule
\end{tabular}
\end{table}

\begin{table}[!h]
\caption{Results on dataset A with 75 customers} \centering
\begin{tabular}{rrrrrrrrr}
\toprule
&\multicolumn{3}{c}{TS}  & \multicolumn{4}{c}{BP}        \\
Instance & NV & TD       & Time(s) & NV & TD       & Time(s) & Gap    \\
\midrule
c101  & 8  & 84.075  & 1.28  & 8  & 84.075  & 1307.64  & 0.00\%  \\
c105  & 8  & 83.800  & 10.61 & 8  & 83.675  & 10707.30 & 0.00\%  \\
c106  & 8  & 80.379  & 0.60  & 8  & 80.743  & 10524.12 & 0.26\%  \\
c107  & 8  & 83.723  & 2.90  & 8  & 83.588  & 10719.11 & 0.06\%  \\
r101  & 16 & 188.867 & 1.50  & 16 & 178.963 & 239.79   & 0.00\%  \\
r102  & 14 & 172.323 & 1.73  & 14 & 164.139 & 481.14   & 0.00\%  \\
r103  & 11 & 146.177 & 19.15 & 11 & 134.750 & 5015.91  & 0.00\%  \\
r105  & 12 & 174.159 & 8.76  & 11 & 154.253 & 744.75   & 0.00\%  \\
r106  & 11 & 146.746 & 8.86  & 10 & 147.330 & 10769.54 & 0.11\%  \\
r109  & 11 & 131.937 & 9.36  & 9  & 135.247 & 6489.55  & 0.00\%  \\
rc101 & 12 & 182.652 & 1.26  & 12 & 176.597 & 2595.03  & 0.00\%  \\
rc102 & 11 & 171.448 & 18.98 & 10 & 168.386 & 10771.39 & 0.10\%  \\
rc103 & 10 & 143.512 & 0.63  & 9  & 146.403 & 10539.07 & 0.84\%  \\
rc106 & 11 & 153.381 & 3.68  & 9  & 158.384 & 10750.79 & 0.20\%  \\
rc107 & 9  & 144.181 & 33.79 & 9  & 144.269 & 10772.05 & 10.87\% \\
rc108 & 9  & 126.413 & 1.85  & 9  & 129.969 & 10144.71 & 12.52\%\\
\bottomrule
\end{tabular}
\end{table}

\newpage
\begin{table}[!h]
\caption{Results on dataset A with 100 customers} \centering
\begin{tabular}{rrrrrrrrr}
\toprule
&\multicolumn{3}{c}{TS}  & \multicolumn{4}{c}{BP}        \\
Instance & NV & TD       & Time(s) & NV & TD       & Time(s) & Gap    \\
\midrule
c101  & 10 & 107.742 & 2.21  & 10 & 107.742 & 3255.73  & 0.00\% \\
r101  & 20 & 211.371 & 6.76  & 19 & 207.751 & 983.72   & 0.00\% \\
r102  & 18 & 206.003 & 18.26 & 17 & 190.955 & 2111.25  & 0.00\% \\
r105  & 15 & 186.240 & 7.66  & 14 & 183.726 & 10374.66 & 5.95\% \\
rc101 & 16 & 225.632 & 5.18  & 14 & 224.172 & 10786.57 & 0.22\% \\
rc102 & 14 & 209.526 & 13.02 & 12 & 221.372 & 10585.25 & 1.80\% \\
rc105 & 15 & 245.247 & 24.93 & 13 & 211.847 & 3935.69  & 0.00\% \\
rc106 & 13 & 195.968 & 22.91 & 12 & 183.270 & 10522.80 & 7.75\%\\
\bottomrule
\end{tabular}
\end{table}
\clearpage \newpage
\section{Impact of Payloads and Slopes on Routing Costs}\label{tab: impact of payloads and slopes full}
\begin{table}[!ht]\caption{Impact of Payloads and Slopes on Routing Costs} 
	\centering 
	  \begin{adjustbox}{max width=\textwidth} 
\begin{tabular}{cc|ccccccccccc|c}
	&         & \multicolumn{12}{c}{Slope factor ($r_2$)}                                                                            \\
	&         & 0.0      & 0.1    & 0.2    & 0.3    & 0.4    & 0.5    & 0.6    & 0.7    & 0.8    & 0.9    & 1.0      & Average \\ \hline
	\multirow{12}{*}{\makecell[l]{Payload\\factor\\ ($r_1$)}} 
	& 0.0       & 1516.0 & 1525.0 & 1543.9 & 1490.0 & 1527.4 & 1466.6 & 1458.7 & 1429.5 & 1456.2 & 1419.1 & 1361.8 & 1472.2  \\
	& 0.1     & 1525.6 & 1550.8 & 1496.3 & 1530.8 & 1536.9 & 1455.9 & 1441.2 & 1480.4 & 1433.8 & 1439.8 & 1375.7 & 1478.8  \\
	& 0.2     & 1514.3 & 1542.1 & 1519.5 & 1536.1 & 1527.3 & 1445.3 & 1472.0 & 1452.6 & 1475.4 & 1445.6 & 1374.8 & 1482.3  \\
	& 0.3     & 1532.6 & 1511.8 & 1564.2 & 1489.7 & 1569.1 & 1478.0 & 1473.6 & 1478.8 & 1470.4 & 1435.7 & 1355.1 & 1487.2  \\
	& 0.4     & 1530.7 & 1545.8 & 1563.0 & 1551.3 & 1540.4 & 1474.5 & 1466.4 & 1454.4 & 1478.4 & 1455.9 & 1363.9 & 1493.2  \\
	& 0.5     & 1530.8 & 1567.3 & 1555.8 & 1567.1 & 1551.9 & 1490.6 & 1470.8 & 1482.0 & 1474.8 & 1439.2 & 1371.1 & 1500.1  \\
	& 0.6     & 1545.9 & 1571.4 & 1565.5 & 1570.8 & 1556.6 & 1470.0 & 1475.9 & 1475.0 & 1512.0 & 1435.6 & 1373.8 & 1504.8  \\
	& 0.7     & 1562.1 & 1541.6 & 1520.5 & 1564.9 & 1534.8 & 1518.6 & 1488.3 & 1489.5 & 1514.7 & 1429.3 & 1377.1 & 1503.8  \\
	& 0.8     & 1553.0 & 1553.9 & 1532.5 & 1564.9 & 1582.0 & 1499.2 & 1511.1 & 1503.4 & 1494.6 & 1465.7 & 1404.7 & 1515.0  \\
	& 0.9     & 1555.6 & 1552.7 & 1568.9 & 1574.3 & 1569.3 & 1498.9 & 1507.3 & 1486.3 & 1495.6 & 1498.7 & 1395.3 & 1518.4  \\
	& 1.0       & 1575.1 & 1569.9 & 1542.8 & 1543.2 & 1586.7 & 1523.5 & 1516.6 & 1544.5 & 1496.1 & 1468.2 & 1413.3 & 1525.4  \\ \hline  
	& Average & 1540.1 & 1548.4 & 1543.0 & 1543.9 & 1553.0 & 1483.7 & 1480.2 & 1479.7 & 1482.0 & 1448.4 & 1378.8 &     
\end{tabular}
\end{adjustbox}
\end{table}

    \end{APPENDICES}

\end{document}